\documentclass[onefignum,onetabnum]{siamonline190516}
\usepackage{amsmath,amssymb,amsfonts,tikz}
\usepackage{url}
\usepackage{etoolbox}
\appto\UrlBreaks{\do\a\do\b\do\c\do\d\do\e\do\f\do\g\do\h\do\i\do\j
\do\k\do\l\do\m\do\n\do\o\do\p\do\q\do\r\do\s\do\t\do\u\do\v\do\w
\do\x\do\y\do\z\do\?\do\1\do\2\do\3\do\4\do\5\do\6\do\7\do\8\do\9\do\0\do\?}
\definecolor{OliveGreen}{cmyk}{0.64,0,0.95,0.60}

\usetikzlibrary{arrows,positioning}

\def\ceil#1{{\lceil #1 \rceil}}
\def\eps{{\varepsilon}}
\newcommand{\R}{\mathbb{R}}
\newcommand{\N}{\mathbb{N}}
\def\e{{\rm e}}

\newtheorem{remark}[theorem]{Remark}

\newenvironment{proofof}{\noindent\sc Proof of}{
    \hspace*{\fill} $\Box$ \vspace{2ex} }

\newcommand{\din}{\delta_{\text{in}}}
\newcommand{\dout}{\delta_{\text{out}}}

\newcommand{\PP}{\textbf{P}}

\newcommand{\argmin}{\operatornamewithlimits{argmin}}

\newcommand{\fin}{f^{\text{in}}}
\newcommand{\fout}{f^{\text{out}}}
\newcommand{\Cin}{C_{\text{in}}}
\newcommand{\Cout}{C_{\text{out}}}
\newcommand{\ain}{\alpha_\text{in}}
\newcommand{\aout}{\alpha_\text{out}}

\newcommand{\beqq}{\begin{equation}}
\newcommand{\eeqq}{\end{equation}}

\ifpdf
\hypersetup{
  pdftitle={On a minimum distance procedure for threshold selection in tail analysis},
  pdfauthor={H. Drees, A. Janssen, S. I. Resnick and T. Wang}
}
\fi

\title{On a minimum distance procedure for threshold selection in tail analysis}

\author{Holger Drees\thanks{Department of Mathematics, University of Hamburg
  (\email{drees@math.uni-hamburg.de}).}
\and Anja Jan\ss en\thanks{Department of Mathematics, KTH Royal Institute of Technology, Stockholm
  (\email{anjaj@kth.se}).}
\and Sidney I. Resnick\thanks{School of Operations Research and
  Information Engineering, Cornell University
  (\email{resnick@cornell.edu}).}
\and Tiandong Wang\thanks{Department of Statistics, Texas A\&M University
(\email{twang@stat.tamu.edu}).}}

\usepackage{amsopn}

\begin{document}

\maketitle

\begin{abstract}%
  Power-law distributions have been widely observed in different areas
  of scientific research. Practical estimation issues include
  selecting a threshold above which
  observations follow a power-law distribution and then  estimating the power-law tail index.  A
  minimum distance selection procedure (MDSP) proposed by Clauset et al. (2009) has been widely adopted in
  practice for  the analyses of social networks.  However, theoretical justifications
  for this selection procedure remain scant.  In this paper, we study
  the asymptotic behavior of the selected threshold and the
  corresponding power-law index given by the
  MDSP. For iid\ observations
    with Pareto-like tails, we derive the limiting distribution of the
    chosen threshold and the power-law index estimator, where the
    latter estimator   is not asymptotically normal. We deduce that in this iid setting
    MDSP tends to choose too high a threshold level and show with asymptotic analysis and simulations how the variance increases compared to Hill estimators
    based on a non-random threshold. We also provide simulation
    results for dependent preferential attachment network data and find that the performance of the
    MDSP procedure is highly dependent on the chosen model
    parameters.
\end{abstract}

\begin{keywords}
Power Laws, Threshold Selection, Hill Estimators, Empirical Processes, Preferential Attachment
\end{keywords}

\begin{AMS}
62G32, 60G70, 62E20, 60G15, 62G30, 05C80
\end{AMS}

\section{Introduction}

In empirical studies, it is common to see observations cluster around
a typical value and use mean
and standard deviation to summarize the distribution. However,
not all distributions satisfy this
pattern and often extreme data not adequately summarized by
  moments is of critical importance.  Power-law distributions give
  emphasis to extreme values in the data and therefore
have attracted scientific interest. These have
been used to model a wide variety of physical, biological, and
man-made phenomena, e.g.\ the sizes
of power outage \cite{CLDN02}, the foraging pattern of various species
\cite{Hetal10}, the frequencies of family names \cite{ZM01}, wealth and income distributions \cite{mandelbrot, yakovenko} and in-
and out-degrees in social networks \cite{BBCR03, BA99, KR01}. A survey
of power laws with focus on network science is given in
\cite{M04}. There is also parallel literature
\cite{Beiretal04,coles:2001,dehaan:ferreira:2006,resnickbook:2007}
 in the extreme value and
heavy-tail community more focussed on environmental science, insurance
and finance.

Power laws are ubiquitous in social network modeling.
 Data repositories of large network data sets such as KONECT \cite{kunegis:2013}, provide estimates of power-law indices
as one of the key summary statistics for almost all listed networks.
\begin{figure}\begin{center}
    \includegraphics[width=\textwidth]{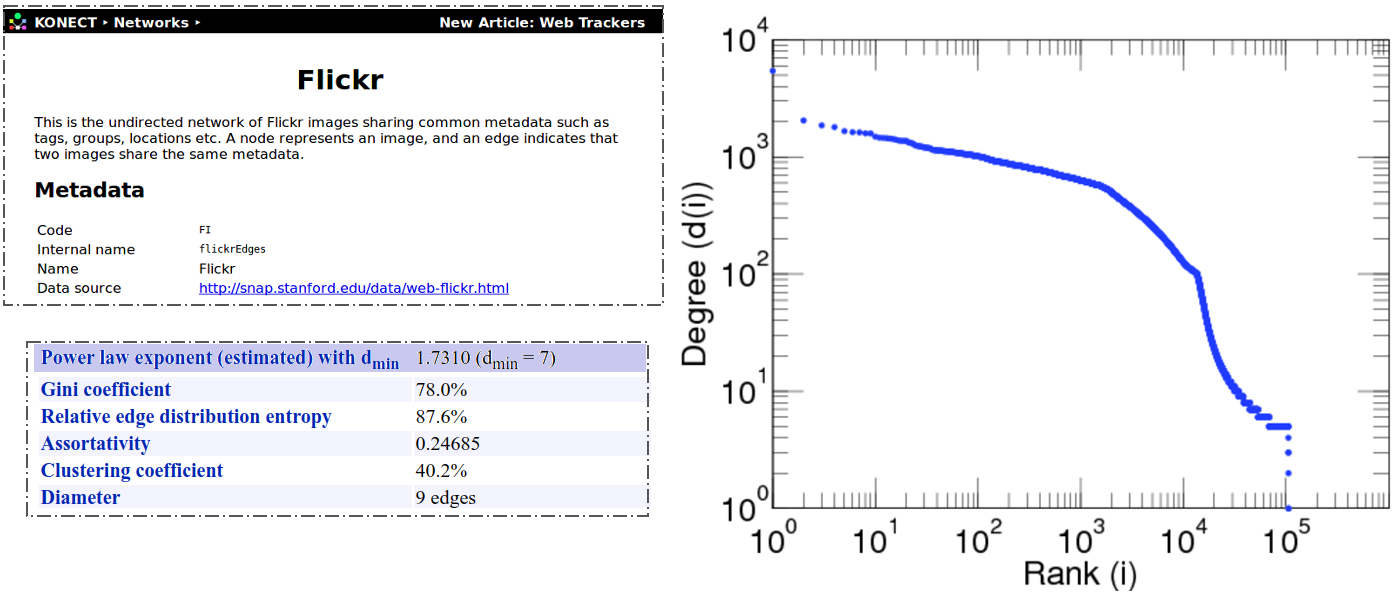}
  \end{center}
  \caption{Snapshots from \url{http://konect.cc/networks/flickrEdges}
    and \url{http://konect.uni-koblenz.de/networks/flickrEdges} for
    the Flickr image metadata data set.}
  \label{fig:flickr}
\end{figure}
\Cref{fig:flickr} displays a snapshot of some of the
  summary statistics of the Flickr image metadata from the KONECT
  database on the left-hand side and a plot of the degree versus the
  rank, both on a logarithmic scale, on the right-hand side. This
  so-called Zipf plot clearly shows that for the 1000 largest degrees
  the logarithm of the inverse of the empirical cdf is well
  approximated by a linear function of the logarithm of its
  argument. In other words, above a threshold of about 500, the tail
  of the cdf  is approximately equal to a multiple of a power
  function, $1-F(x)\approx  cx^{-\alpha}$, $x\ge 500$, where $\alpha$
  can be estimated by the reciprocal of the
  absolute value of the slope of a line fitted
  to the left part of the plot.

  Power-law index estimation is an important task for analysis of many
  heavy tailed phenomena as the power-law index controls
  important characteristics. In social networks,
  the power-law degree index controls the likelihood of a node having
  large degree. It also describes the growth rate of
  degrees and maximal degree
  as a function of the number of edges.

In social network analyses, one popular modeling choice is the
preferential attachment (PA) model \cite{BBCR03, durrett:2010, KR01,
  krapivsky:2001, vanderHofstad:2017}, where nodes and edges are added
according to the probabilistic rule that nodes with larger degrees
tend to attract more edges. Limit theory for degree counts in a linear PA
model can be found in \cite{bhamidi:2007, KR01, krapivsky:2001,
  resnick:samorodnitsky:2015, resnick:samorodnitsky:2016b,
  resnick:samorodnitsky:towsley:davis:willis:wan:2016,
  wang:resnick:2016, wang:resnick:2015, wang:resnick:2018a,
  wang:resnick:2018b} where in particular it is shown that the linear PA network
models generate power-law degree distributions in the limit. The corresponding tail indices
are related to the offset parameter, and with power-law indices of the in- and out-degree distribution estimated,
one can infer values of the other parameters. In \cite{wan:wang:davis:resnick:2019} it is shown in simulations that this approach using extreme value techniques is more robust than a parametric maximum likelihood estimator against modeling errors.

When one has identified a threshold above which the empirical degree
distribution can be well approximated by a power law, the tail index
can  be estimated using the  Hill estimator; see
subsection \ref{subsec:background} for a brief introduction. However,
the performance of this estimator strongly depends on the chosen
threshold. On the one hand, if it is chosen too large, then the
estimator uses too little data which results in a high variance. On
the other hand, if one selects too small a threshold, say 200 for the
Flickr friendship data, then one tries to fit a line to a non-linear
relationship between the rank and the degree in the Zipf plot, which
will lead to a large bias. Identifying a range on which
one may fit the plot by a straight line is an important step to
understand the tail behavior of the degree distribution.

A threshold selection procedure is proposed in \cite{CSN09}, where one
chooses the cutoff value that yields the smallest
Kolmogorov-Smirnov distance between the empirical distribution above
the threshold and the corresponding fitted power law.
This selection procedure combined with the Hill estimator is for
instance used by KONECT \cite[page 30]{kunegis:2018}
to compute the power-law exponent shown in
\Cref{fig:flickr}, with $d_{min}$ denoting the selected threshold. It
has been widely adopted in the analyses of social networks and income
distributions, having attracted more than 3,000 citations;
see, for example, \cite{ahn2007analysis, cho2011friendship, java2007we, kivela2014multilayer, leskovec:2009,OANCEA2017486, SAFARI2018169, SORIANOHERNANDEZ2017403}. It has also been encoded as an R-package called \verb6poweRlaw6 (cf. \cite{gillespie:2015}). We now outline this threshold selection method.

\subsection{Minimum distance selection procedure (MDSP)}\label{subsec:background}
Mathematically, a nonnegative random variable $X$  follows a
power-law distribution if
its tail distribution function satisfies
\begin{equation} \label{eq:Partail}
 1-F(x)  = cx^{-\alpha},
\end{equation}
for $x$ exceeding some threshold $x_0 > 0$, where $c > 0$ is some constant and $\alpha > 0$ is known as the exponent or tail index. The distribution that fulfills this relation for all $x>c^{1/\alpha}$ is called Pareto distribution. Therefore, it is also common to speak of a Pareto tail instead of a power tail.

As mentioned in \cite{CSN09}, empirical
distributions rarely  follow a power law for
all values, but rather only for observations greater than some cutoff
value. Therefore, there are
two parameters to determine: the exponent $\alpha$ and the cutoff value $x_0$. Provided that we have a
good estimate for the threshold $x_0$, we can discard all observations below $x_0$ and estimate $\alpha$ by the
maximum likelihood estimator based on the remaining exceedances. For iid\ observations with distribution function \eqref{eq:Partail}, the exceedances over $x_0$, divided by $x_0$, follow a Pareto distribution with parameters $\alpha $ and $c=1$. Taking the logarithm of those observations results in a sample of iid\ exponential random variables with mean $1/\alpha$, leading to a straightforward maximum likelihood estimator of $\alpha$. As the true value of $x_0$ is unknown, it has to be chosen depending on the data, and often a suitable order statistic of the observations is used. More precisely, suppose random variables  $X_i$, $1\le i\le n$, are
observed and denote the order statistics by $X_{1:n}\le
X_{2:n}\le\cdots\le X_{n:n}$. If one uses the $k$th largest order
statistic as a threshold, the maximum likelihood approach applied to
model \eqref{eq:Partail} leads
to the well-known Hill estimator \cite{Hill75}
\begin{equation} \label{eq:Hilldef}
   \hat\alpha_{n,k} := \left( \frac 1{k-1} \sum_{i=1}^{k-1} \log \frac{X_{n-i+1:n}}{X_{n-k+1:n}}\right)^{-1}.
\end{equation}
For iid\ data which is generated from a distribution of the
  form \eqref{eq:Partail}, this estimator is known to be consistent
\cite{Mason82} and asymptotically normal \cite{Hall82} with rate
$k^{-1/2}$, provided $X_{n-k:n}$ exceeds the threshold with
probability tending to 1. Consistency and normality also
  hold in the more general case of iid\ observations coming from a
  regularly varying distribution under suitable assumptions about the
  growth of $k$ and control of bias. For discussion,
    see
  \cite{dehaan:ferreira:2006,resnickbook:2007}. Furthermore, results
can be extended to the case of dependent observations from stationary
sequences under suitable mixing conditions; see e.g.\ \cite{Hsing91,Drees2000}. For the case of
dependent data consisting of node degrees from
  preferential attachment models, the consistency of the Hill
estimator has recently been shown in
\cite{wang:resnick:2018a,wang:resnick:2018b}
  but asymptotic normality of the Hill
  estimator in this setting remains an open problem. In any case, for the iid, stationary mixing and
  network settings, performance of the Hill estimator critically
  depends on $k$, that also determines the estimator $X_{n-k+1:n}$ of $x_0$.

Clauset et al. \cite{CSN09} suggest estimating the
cutoff value by the order statistic which minimizes the Kolmogorov-Smirnov
distance between the empirical distribution of the exceedances and the
Pareto distribution fitted with the larger order statistics. To be more
precise, define the Kolmogorov-Smirnov distance
\begin{equation} \label{eq:KSdist}
 D_k := \sup_{y\ge 1} \Big| \frac 1{k-1} \sum_{i=1}^n 1_{(y,\infty)}\Big( \frac{X_{n-i+1:n}}{X_{n-k+1:n}}\Big)-y^{-\hat\alpha_{n,k}}\Big|
\end{equation}
and use $X_{n-k_n^*+1:n}$
with
\begin{equation} \label{eq:KSkopt}
  k_n^* := \argmin_{k\in\{2,\ldots,n\}} D_k
\end{equation}
as an estimator of the unknown threshold. (If the point of minimum is
not unique, we may e.g.\ choose the smallest one.)
Since we choose the threshold that minimizes the
 distance between fitted and empirical tail this method is called the
 minimum distance selection procedure (MDSP).  This
  method has also been adapted to binned data in \cite{VC14}.

\subsection{Summary}
\label{sec:sum}

The MDSP method is widely applied in practice, particularly in
  Computer Science and Network Science. However, to the
best of our knowledge its performance has not been mathematically
analyzed even in classical contexts where data is assumed to come from
an iid model of repeated sampling. Under the assumption of iid
observations, we will show:
\begin{itemize}
  \item The MDSP tail index estimator has a
    non-normal limit distribution that is difficult to calculate
    rendering the task of computing reliable confidence intervals
    quite difficult.
\item The MDSP procedure often leads to
choosing a $k^*_n$ that is too small, resulting in increased
  variance and root mean squared error (RMSE) for the Hill estimator relative
to a choice which minimizes the asymptotic mean squared
error.
\end{itemize}
For dependent observations simulated from a preferential
  attachment model, we observe that the performance of the MDSP
  estimator strongly depends on the chosen
  model parameters, with certain  sets of model
    parameters leading to good MDSP performance but other choices of
    model parameters leading to simulations where the MDSP choice of
    $k_n^*$ is too large. Mathematical analysis of MDSP
      performance for preferential attachment data is not available
      since the Brownian motion embedding techniques successful for
      the iid case are not available for the network case.

We begin in
Section \ref{sec:Pareto}  with
the iid case assuming the underlying distribution is exact Pareto and
thus $k=n$ would be the best choice for minimizing asymptotic RMSE that is associated with Cram\'er-Rao    bounds and MLE asymptotic efficiency.
It will be shown that the
distribution of $k_n^*/n$ can be approximated by a distribution
supported by the whole interval $(0,1]$, so that with non-negligible
probability $k_n^*$ is much smaller than $n$. In Section 3, we analyze
the asymptotic behavior of $k_n^*$ if the underlying cdf satisfies
\eqref{eq:Partail} for all $x>x_0$ for some $x_0$ such that
$F(x_0)>0$, but $F$ shows a different behavior below $x_0$. In Section
\ref{sec:networks}, we discuss numerical results for the
performance of the MDSP applied to
the in-degrees of linear preferential attachment
networks. All proofs are postponed to the appendix.

The MDSP offers attractive features.  The procedure
yields estimates without requiring user discretion. It is readily
implemented with R-packages that are well designed and can be ported into
another algorithm. In network simulations the tail index estimates provided by MDSP have often proved to be
reasonable, provided network parameters are close to those observed in
empirical studies, cf.\ Section \ref{subsec:sim}. However, this method has limitations and needs to
be applied with caution. Even in the  classical iid case,
the asymptotic theory of MDSP estimation is fairly complex and it is
not an easy task to  extract confidence intervals for estimates
obtained by this method. Furthermore MDSP estimates of the tail index do not achieve minimal
asymptotic RMSE. For the node based data of random graphs, there is no
theoretical analysis available for MDSP estimates.

\section{The Pareto case}
\label{sec:Pareto}

Throughout this section, we assume that the observations are independently  drawn from an exact Pareto distribution on the interval $(c^{1/\alpha},\infty)$, that is
\begin{equation}  \label{eq:Parcdf}
  1-F(x) = cx^{-\alpha}, \quad x\ge c^{1/\alpha},
\end{equation}
for some $\alpha,c>0$. Such a model rarely arises in practice, but we will see that one of
the main drawbacks of the MDSP can most easily be explained in this
setting. Moreover, in Section \ref{sec:Parbreak}, this drawback will
be observed in a modified form in more complex and realistic
models. In the case of \eqref{eq:Parcdf}, the sum in
  \eqref{eq:Hilldef} consists of $k-1$ summands which are all
  approximately iid\ exponential random variables with mean
  $1/\alpha$. A reasonable selection procedure should therefore use as many observations as possible in order to minimize the mean squared error of the Hill estimator. This is obvious from the fact that the rescaled exceedances $X_{n-i+1:n}/X_{n-k+1:n}$, $1\le i\le k-1$, have the same distribution as the order statistics of $k$ Pareto rv's with $c=1$, so that choosing $k<n$ is equivalent to reducing the sample size.

Note that
\begin{align} \label{eq:Dk}
  D_k & = \sup_{y\ge 1} \bigg| \frac 1{k-1} \sum_{i=1}^{k-1} 1_{(y,\infty)} \Big(\frac{X_{n-i+1:n}}{X_{n-k+1:n}}\Big)-y^{-\hat\alpha_{n,k}}\bigg|\nonumber\\
  & = \max_{1\le j< k}\max\bigg( \Big(\frac{X_{n-j+1:n}}{X_{n-k+1:n}}\Big)^{-\hat\alpha_{n,k}}-\frac{j-1}{k-1},
  \frac j{k-1}-\Big(\frac{X_{n-j+1:n}}{X_{n-k+1:n}}\Big)^{-\hat\alpha_{n,k}}\bigg)\nonumber\\
  & = \max_{1\le j\le k} \bigg|\Big(\frac{X_{n-j+1:n}}{X_{n-k+1:n}}\Big)^{-\hat\alpha_{n,k}}-\frac{j}k\bigg|+rem
\end{align}
with $rem$ denoting a remainder term with modulus of at most $1/k$.
It is well known that $n^{1/2}D_n$ weakly converges to the supremum of a Brownian bridge if the Hill estimator is replaced with the true value $\alpha$. More generally, Theorem \ref{th:Parmain} given below shows that $n^{1/2} D_{\ceil{nt}}$ converges to $\sup_{s\in(0,1]}|Z(s,t)|$ for some Gaussian process uniformly for all $t\in[\eps,1]$ for any $\eps\in(0,1)$. The limit process is self-similar with $\sup_{s\in(0,1]}|Z(s,t)|=^d t^{-1/2} \sup_{s\in(0,1]}|Z(s,1)|$. Hence, it is more likely that its infimum is attained at some $t$ close to 1 than in the neighborhood of some smaller value.  However, with non-negligible probability the point of minimum of $t\mapsto\sup_{s\in(0,1]}|Z(s,t)|$ is considerably smaller than 1, corresponding to a suboptimal behavior of the MDSP.

\begin{theorem} \label{th:Parmain}
 If $X_i$, $i \in \mathbb{N}$, are iid with cdf $F$ given in \eqref{eq:Parcdf}, then $D_k$ as defined in \eqref{eq:KSdist} satisfies (for suitable versions of $X_i$)
  \begin{equation} \label{eq:Dntapprox}
     n^{1/2} D_{\ceil{nt}}  = \sup_{0<s\le 1} |Z_n(s,t)|+ O_P\Big(\frac{\log(nt)(\log(nt)+(\log n)^{1/2})}{n^{1/2}t}\Big)
  \end{equation}
  uniformly for $t\in [2/n,1]$ with
  $$
    Z_n(s,t)  :=  \Big(\frac{W_n(st)}{st}-\frac{W_n(t)}t\Big)s  + \bigg(\int_0^1 \frac{W_n(tx)}{tx}\, dx -\frac{W_n(t)}t\bigg)s\log s,
  $$
  where $W_n, n \in \mathbb{N},$ is a suitable sequence of Brownian motions.
\end{theorem}

\begin{figure}[tb]
  \centering

  \includegraphics[width=13cm]{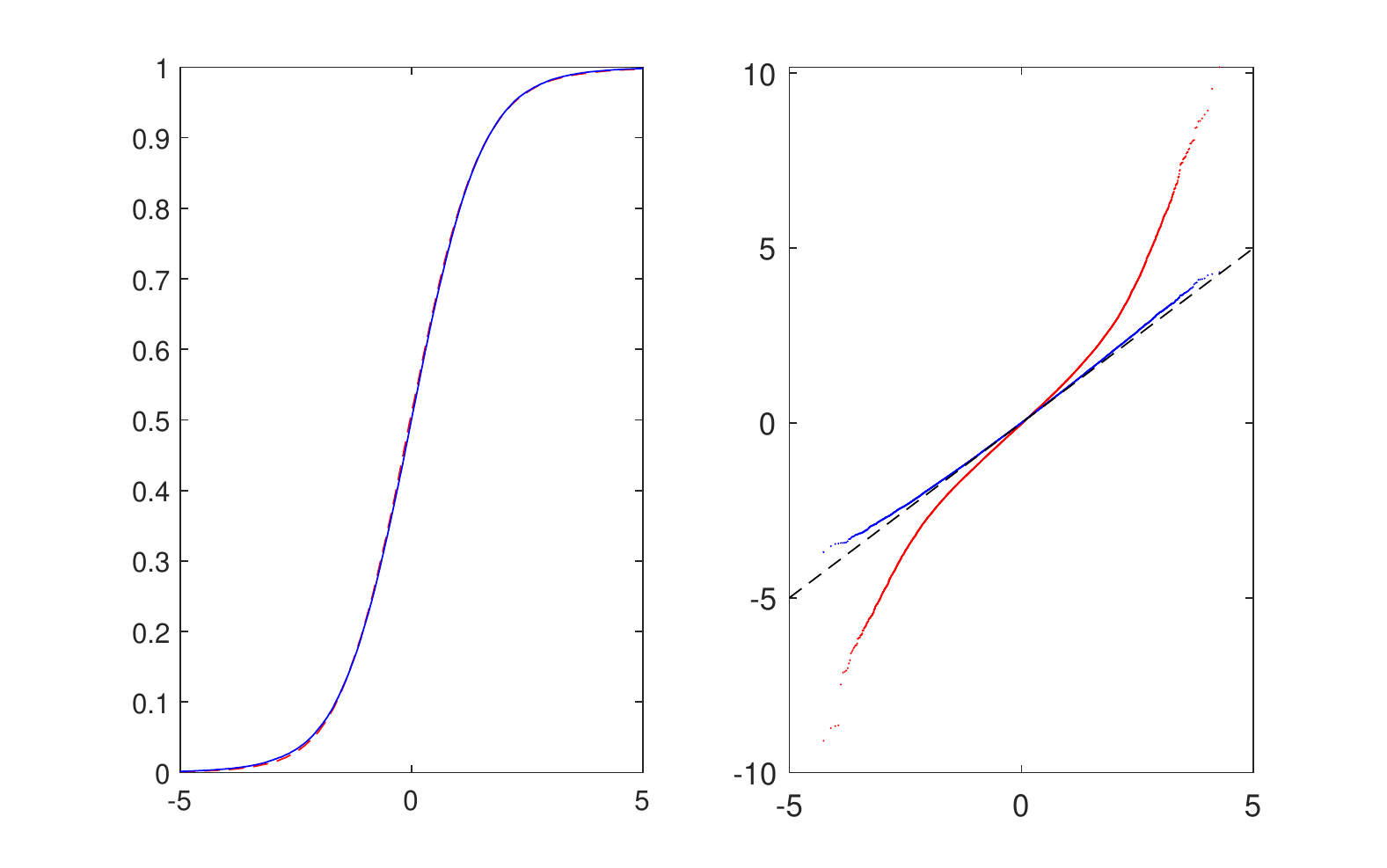}
  \caption{Left: empirical cdf of $n^{1/2}(\hat\alpha_{n,k_n^*}-\alpha)$ for $n=1,000$ (red, dashed) and limit cdf according to \eqref{eq:Parasympt} (blue, solid); right: normal Q-Q plot of $n^{1/2}(\hat\alpha_{n,k_n^*}-\alpha)$ (red) and of $n^{1/2}(\hat\alpha_{n,n}-\alpha)$ (blue) for $n=1,000$, the black dashed line is the main diagonal.}
  \label{fig:HillcompPar}
\end{figure}

From Theorem \ref{th:Parmain}
we obtain the joint asymptotic distribution of the selected number $k$
and the resulting Hill estimator provided the limiting process has a unique
point of minimum:
\begin{corollary} \label{cor:Par}
  If $t\mapsto\sup_{0<s\le 1}|Z_1(s,t)|$ has a unique point of minimum
  $T$ a.s.,
  then with $k_n^\ast$ as in \eqref{eq:KSkopt} and $\hat{\alpha}_{n,k}$ as in \eqref{eq:Hilldef} we have
  \begin{equation} \label{eq:Parasympt}
   \big( k_n^*/n,n^{1/2}(\hat\alpha_{n,k_n^*}-\alpha)\big) \,\longrightarrow  \bigg(T, \alpha\Big(\int_0^1 \frac{W_1(Tx)}{Tx}\, dx -\frac{W_1(T)}T\Big)\bigg) \quad \text{weakly}
  \end{equation}
  as $n\to\infty$.
\end{corollary}
Note that while the second component of the limit would be normally distributed if $T$ were deterministic, due to the randomness of $T$ the MDSP tail index estimator $\hat\alpha_{n,k_n^*}$ does not have a normal limiting distribution; see \Cref{fig:HillcompPar} for a plot of its limiting cdf and a Q-Q plot that clearly shows the strong deviation from normality.

\begin{remark}{\rm
  Unfortunately, the standard techniques to prove the uniqueness of the point of minimum of a Gaussian process apparently do not carry over to our limit process. However, the simulations outlined below suggest that indeed \eqref{eq:Parasympt} holds. In any case, Theorem 3 of \cite{Ferger04} implies the following weaker result:

  If all points of minimum of $t\mapsto\sup_{0<s\le 1}|Z_1(s,t)|$ lie in a random interval $[T_0,T_1]$, then
  $$ P\{T_1<x\} \le \liminf_{n\to\infty} P\{\hat k_n^*/n<x\} \le \limsup_{n\to\infty} P\{\hat k_n^*/n\le x\}\le P\{T_0\le x\}
  $$
  for all $x\in[0,1]$.}
\end{remark}

\begin{figure}[tb]
  \centering

  \includegraphics[width=10cm]{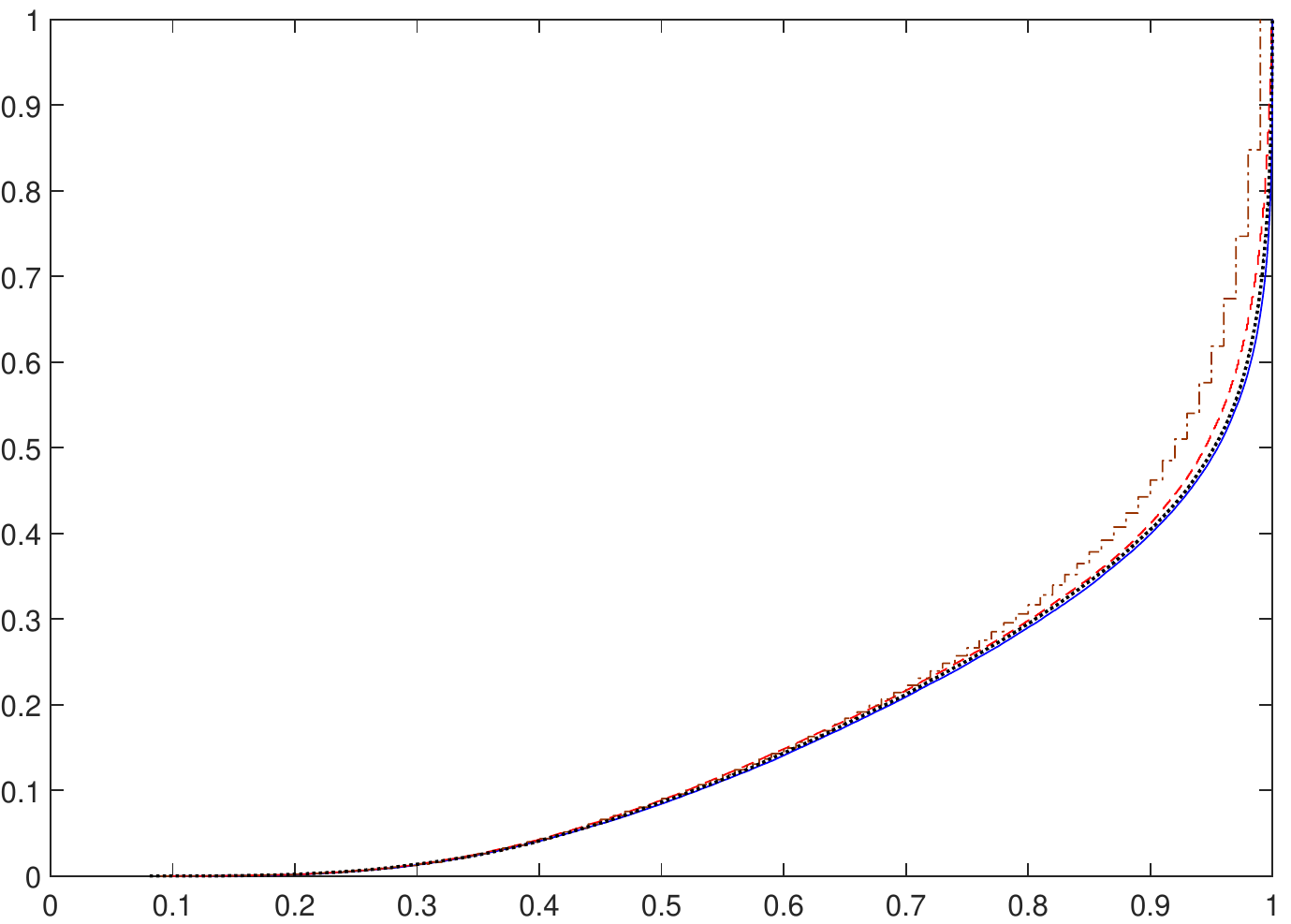}
  \caption{Empirical cdf of $k_n^*/n$ for $n=100$ (brown, dash-dotted), $n=1,000$ (red, dashed) and $n=10,000$ (black, dotted) and limit cdf according to \eqref{eq:Parasympt} (blue, solid) for a Pareto model with $\alpha=c=1$}
  \label{fig:optkcomp}
\end{figure}

\Cref{fig:optkcomp} shows the empirical cdf of $k_n^*/n$ calculated from $10^5$ simulations of standard Pareto samples (i.e., $\alpha=c=1$) of size $n\in\{100; 1,000; 10,000\}$ in comparison with the limit cdf from \eqref{eq:Parasympt}. (The latter was approximately calculated from $10^5$ simulations of a discretized version of the limit process $Z_1$ on a grid with $5\cdot 10^4$ points in each argument.) The difference between the cdf of $k^*_n/n$ and the limit cdf of $T$ is small for $n=1,000$ and hardly visible for $n=10,000$, while $k_n^*/n$ is stochastically a bit smaller for $n=100$.

Reading from Figure \ref{fig:optkcomp},  in the limit, the
probability that $k_n^*$ is less than $(3/4)n$ is about $1/4$, and the
corresponding probabilities for $n/2$ and $n/3$ are $8.4\%$ and $2\%$,
respectively. So while in about $3/5$ of all cases $k_n^*$ is at least
$0.9n$, there is a non-negligible probability that $k_n^*$ is
substantially smaller than $n$. As a consequence, the variance and the
root mean squared error (RMSE) of the corresponding Hill estimator
$\hat\alpha_{n,k_n^*}$ are much larger than those of the
estimator $\hat\alpha_{n,n}$ with minimal RMSE. In the limit, the variance of
$\hat\alpha_{n,k_n^*}$ is about $88\%$ larger, resulting in an RMSE
which is $37\%$ higher. For finite sample sizes the corresponding
figures for the RMSE are $33\%$ for $n=100$, $36\%$ for $n=1,000$ and $37\%$ for
$n=10,000$.

The left plot of \Cref{fig:HillcompPar} compares the (empirical)
distribution of $n^{1/2}(\hat\alpha_{n,k_n^*}-\alpha)$ for $n=1,000$
with the limit distribution given in \eqref{eq:Parasympt}. Here the
approximation is even better than in Figure \ref{fig:optkcomp}.
The right
    plot shows a normal Q-Q plot of the standardized estimation errors
    of $\hat\alpha_{n,k_n^*}$ and $\hat\alpha_{n,n}$,
    respectively. While the latter estimator is asymptotically normal
    (as is $\hat\alpha_{n,\ceil{nt}}$ for all $t\in(0,1]$), the Hill
    estimator based on the top $k_n^*$ order statistics has much
    heavier tails.

The heavier tail of $\hat \alpha_{n,k^*_n}$ is important when
constructing confidence intervals.
The analyses in \cite{CSN09} tempt to use $I_{1-\beta} :=
[\hat\alpha_{n,k_n^*}(1-(k_n^*)^{-1/2}c_{\beta/2});
\hat\alpha_{n,k_n^*}(1+(k_n^*)^{-1/2}c_{\beta/2})]$, with $c_{\beta/2}=\Phi^{-1}(1-\beta/2)$
being the standard normal quantile to the level $1-\beta/2$, as
the confidence interval with
asymptotic level $1-\beta$.
This approach, however,
ignores the inherently stochastic nature of $k_n^*$ and leads to a
severe underestimation of the actual error. For example, the
non-coverage probability of $I_{0.95}$ is greater than $8\%$ and the
one of $I_{0.99}$ larger than $2\%$ for $n\in\{100;1,000\}$. In
contrast, the confidence interval derived from \eqref{eq:Parasympt}
$$ I_{1-\beta}^* := [\hat\alpha_{n,k_n^*}(1-n^{-1/2}c_{\beta/2}^*); \hat\alpha_{n,k_n^*}(1+n^{-1/2}c_{\beta/2}^*)]
$$
with $c_{\beta/2}^*$ denoting the $(1-\beta/2)$-quantile of the limit distribution of the standardized estimation error is very accurate. For $\beta=5\%$ and $c_{\beta/2}^*\approx 2.74$ the non-coverage probability is $4.8\%$, for $\beta=1\%$ and $c_{\beta/2}^*\approx 4.09$ it equals $1\%$.

\section{Deviations from the Pareto model}
\label{sec:Parbreak}

In \cite{CSN09}, the authors assess the accuracy of the MDSP by examining whether the
method is able to recover the threshold level. One example used in \cite{CSN09} is the simulated data which
follow a pure Pareto distribution beyond the threshold but exponential below.
We now provide more general analyses and consider different cases in which the data distribution deviates
from the pure Pareto.

Throughout this section, we assume that iid rv's $X_1, \ldots, X_n$ are observed with a cdf of the type \eqref{eq:Parcdf} above a threshold, whereas below the cdf deviates from this power function. For the asymptotic analysis, it is most convenient to specify the model in terms of the quantile function $F^\leftarrow (1-t)$, $0<t<1$. In the pure Pareto case analyzed in the previous section, \eqref{eq:Parcdf} leads to $F^\leftarrow(1-t) = (t/c)^{-1/\alpha}$, $0<t<1$. In the following we will analyze the case in which the quantile function has this form only above a certain threshold. To this end, let
\begin{equation} \label{eq:Parbreakcdf}
  F^\leftarrow(1-t) = (t/c)^{-1/\alpha} (1+H(t)) \quad \text{with}
  \quad H(t)=0\;\forall\, t<t_0,
\end{equation}
for some $0<t_0<1$. If $H(t)$ is different from 0 in a right neighbourhood of $t_0$ this implies a deviation from the pure Pareto model. We will investigate whether the minimum distance selection procedure will be able to detect the change point $t_0$ accurately by estimating a value of $k_n^\ast$ such that $k_n^\ast/n$ is close to $t_0$.

The smoothness of the function $H$ at $t_0$ is crucial for the tail behavior. If $H$ jumps at $t_0$, then with growing sample size the resulting discontinuity of $F^\leftarrow$ at $1-t_0$ could easily be identified from the empirical distribution function, making it thus even possible to spot the break point $t_0$ with the naked eye. Therefore, we focus on the more challenging case that $H$ is continuous at $t_0$. On the other hand, if $H$ is very smooth in a neighborhood of $t_0$ then the behavior of the quantile function hardly changes at $t_0$ and no clear break point exists. In that case, one cannot expect the MDSP to find a value of $k_n^\ast$ such that $k_n^\ast/n$ is close to $t_0$. Moreover, in this setting the exact value of $t_0$ would be of less practical relevance, since the exceedances would not significantly deviate from a Pareto distribution even if one chooses a substantially lower threshold. For that reason, we will in Theorem \ref{th:Parbreak} analyze the asymptotic behavior of $k_n^*/n$ under the assumption that $H$ has a non-vanishing right hand derivative at $t_0$, ensuring therefore a smooth but detectable change in the behavior of the tail at $t_0$. Indeed, the Zipf plot in \Cref{fig:flickr} suggests such a continuous quantile function with a kink at some point slightly above rank 1000. In addition to our theoretical analysis of this case, the performance of the MDSP will be examined for examples of the aforementioned cases of discontinuous or very smooth functions $H$ in a small simulation study.

Note that the number of observations above the threshold
$F^\leftarrow(1-t_0)=(t_0/c)^{-1/\alpha}$ is binomially distributed
with parameters $n$ and $t_0$ and has mean $nt_0$ and these are
  the observations that come from a distribution with exact Pareto tail. Hence, by the normal approximation to the binomial
  distribution, for any sequence $l_n=o(n)$
such that $n^{1/2}=o(l_n)$, with probability tending to 1 the
Kolmogorov-Smirnov distance $D_k$ will be exactly the same as in the
Pareto case for all $k\le nt_0-l_n$. Conversely, if $k>nt_0+l_n$, then
a substantial part of the observations will be smaller than the
threshold and,  due to the deviation
from the Pareto model, the distance $D_k$ will be large. (In fact, this is the main motivation for the
MDSP.) Therefore, we expect that with probability tending to 1, the minimal value of
$D_k$ will not be assumed on the range $\{nt_0+l_n,\ldots,n\}$. So the
most interesting case arises when $k$ deviates from $nt_0$ by the
order $n^{1/2}$. We will see that then a new effect emerges: the
influence of the deviation from the Pareto model on $D_k$ may partly
cancel out the random deviation of the fitted Pareto cdf from the
empirical cdf (in the pure Pareto model), leading to an overall
smaller Kolmogorov-Smirnov distance. As this effect can only occur if
$k$ is sufficiently close to $nt_0$, it increases the probability that
$k_n^*$ takes on some value in that range. Nevertheless, as in the
Pareto case, quite often the MDSP picks too small a value which leads
to a substantial increase in the RMSE of the resulting Hill estimator compared with a choice
of $k$ minimizing the RMSE.

The following analog to Theorem \ref{th:Parmain} is
weaker than what could be stated but avoids some technicalities.
\begin{theorem} \label{th:Parbreak}
  Suppose \eqref{eq:Parbreakcdf} holds for some function $H$ which is continuous on $(0,1)$ and continuously differentiable on $(t_0,t_0+\delta]$ for some $\delta>0$ with $\lim_{t\downarrow t_0} H'(t)=:h_0\ne 0$. Let $\eps_n\downarrow 0$ be such that $n^{1/2}\eps_n\to \infty$.
  \begin{enumerate}
    \item Approximation \eqref{eq:Dntapprox} holds uniformly for $t\in[2/n,t_0-\eps_n]$.
    \item $n^{1/2} \inf_{k\in\{\ceil{n(t_0+\eps_n)},\ldots,n\}} D_k \to \infty$ in probability.
    \item For all $C>0$,
    \begin{equation} \label{eq:Dntapprox2}
     n^{1/2} D_{\ceil{nt_0+n^{1/2}u}}  = \sup_{0<s\le 1} |\tilde Z_n(s,u)|+ o_P(1)
  \end{equation}
  uniformly for $u\in [-C,C]$ with
  $$
    \tilde Z_n(s,u)  :=  Z_n(s,t_0) +\alpha h_0\big(u+W_n(t_0)-t_0W_n(1)\big)^+s(1+ \log s).
  $$
  \end{enumerate}
\end{theorem}
So while the limit process is unchanged for $t<t_0$, for $k\sim nt_0$ the behavior of $k$ on the finer scale $n^{1/2}$ (instead of $n$) influences the asymptotic behavior of $D_k$. Since $u$ can be arbitrary, one  expects that the asymptotic behavior of $k_n^*/n$ changes in comparison with the Pareto case discussed in Section \ref{sec:Pareto} as follows.

Let
$$ z_{t_0,\min} := \inf_{u\in\R} \sup_{0<s\le 1} |\tilde Z_1(s,u)| $$
and
$$ T := \argmin_{t\in(0,t_0]} \sup_{0<s\le 1}|Z_1(s,t)|. $$
If $z_{t_0,\min}>\sup_{0<s\le 1}|Z_1(s,T)|$, then $k_n^*\sim nT$, else $k_n^*\sim nt_0$. In the former case, we define $T^*:=T$, in the latter $T^*:=t_0$. Note that $u\mapsto \tilde Z_1(s,u)$ is constant for $u\le t_0W_1(1)-W_1(t_0)$, so that in general there is no unique point of minimum of $\sup_{0<s\le 1} |\tilde Z_1(s,u)|$. Therefore,  Theorem \ref{th:Parbreak} (iii) does not directly
permit conclusions about the fine scale behavior of $k_n^*/n$ in the neighborhood of $t_0$.

However, our simulations suggest that any point of minimum $u$ of $\sup_{0<s\le 1} |\tilde Z_1(s,u)|$  leads to the same value $V^*:=\alpha h_0(u+W_1(t_0)-t_0 W_1(1))^+$. The proof of Theorem \ref{th:Parbreak} shows that the limit distribution  of the Hill estimator $\hat\alpha_{n,k_n^*}$ only depends on $T^*$ and $V^*$; cf.\ \eqref{eq:alphaapprox} and \eqref{eq:apprfactor2}. For that reason, we conjecture that
\begin{equation} \label{eq:HillParbreakasymp}
  \big( k_n^*/n,n^{1/2}(\hat\alpha_{n,k_n^*}-\alpha)\big) \,\longrightarrow  \bigg(T^*, \alpha\Big(\int_0^1 \frac{W_1(T^*x)}{T^*x}\, dx -\frac{W_1(T^*)}{T^*}+V^*1_{\{T^*=t_0\}}\Big)\bigg)
\end{equation}
weakly.

\begin{figure}[tb]
  \centering

  \includegraphics[width=10cm]{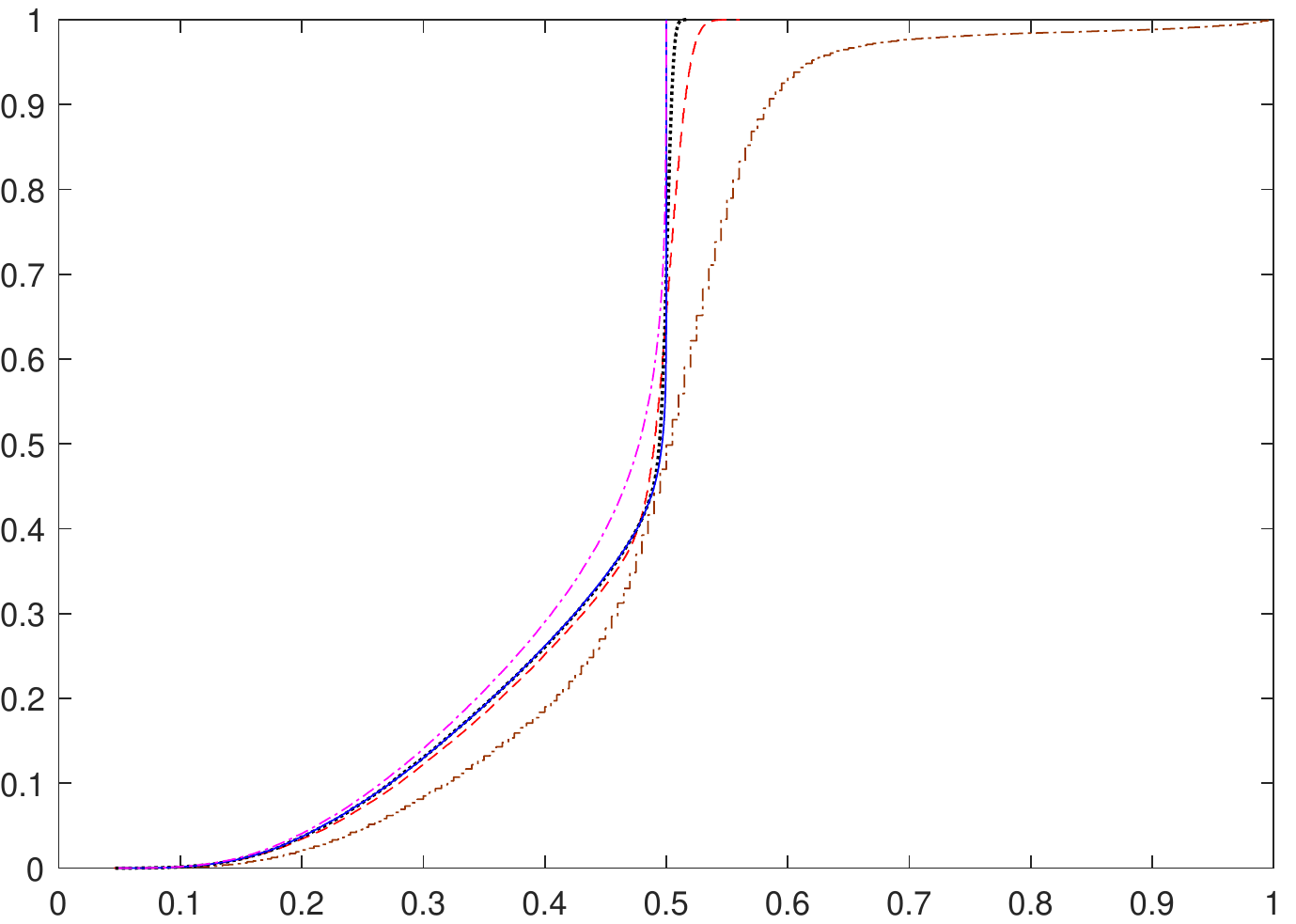}
  \caption{Empirical cdf of $k_n^*/n$ for $n=200$ (brown, dash-dotted), $n=2,000$ (red, dashed) and $n=20,000$ (black, dotted) and limit cdf of $T^*$ (blue, solid) for the piecewise Pareto model described in the text; the magenta dash-dotted line indicates the cdf of $T$, i.e., the properly rescaled limit cdf in the pure Pareto model.}
  \label{fig:optkcompbreak}
\end{figure}

Note that asymptotically $D_k$ shows a different behavior for $k\sim nt$ with $t<t_0$ and for $|k-nt_0|=O(n^{1/2})$. However, for given sample size $n$ and $k<nt_0$, it is not obvious whether to apply Theorem \ref{th:Parbreak} (i) with $t=k/n$ or (iii) with $u=n^{-1/2}(k-nt_0)$. One might expect that this ambiguity is reflected in a loss of accuracy of the approximation of the cdf of $k_n^*/n$ by the cdf of $T^*$, in particular in the vicinity of $t_0$. Figure \ref{fig:optkcompbreak} is the analog to Figure \ref{fig:optkcomp} in the present setting. More precisely, we have chosen $H(t)=\big((t/t_0)^{1/\alpha-1/\beta}-1\big)1_{[t_0,1)}(t)$ such that below the threshold the cdf equals a Pareto cdf with parameter $\beta$ instead of $\alpha$. The parameters are chosen as $\alpha=1$, $\beta=1/2$ and $t_0=1/2$. To make the simulations comparable, the sample sizes are now chosen to be $n\in\{200;2,000;20,000\}$ such that the expected number of observations drawn from the Pareto tail is the same as in Figure \ref{fig:optkcomp}. In addition, the magenta dotted line shows the cdf of $T$, i.e.\ the limit cdf in the Pareto case rescaled to the interval $[0,t_0]$. As explained above, the new limit distribution of $T^*$ has positive mass at $t_0$, leading to a substantial shift of the distribution towards $t_0$ and thus an improved performance of the MDSP.

The distribution of $k_n^*/n$, being pre-asymptotic, smears this mass over a neighborhood of $t_0$; this effect is the stronger the smaller $n$ is (as one has expected since the neighborhood is of the order $n^{-1/2}$). For $n=20,000$, again the approximation of the cdf of $k_n^*/n$ by the limit is almost perfect everywhere; for $n=2,000$ it is very good up to the $70\%$-quantile while (unlike the limit) the cdf of $k_n^*/n$ has  some mass above $t_0$, resulting in a visible approximation error in the upper tail. For $n=200$, the approximation is overall quite  poor. Observe that while in the Pareto case the cdf of $k_n^*/n$ is stochastically increasing with the sample size, here it is stochastically decreasing with $n$, because for smaller sample sizes the Kolmogorov-Smirnov test detects the deviation from the Pareto tail later and this effect is much more pronounced than the differences in the Pareto case.

\begin{figure}[tb]
  \centering

  \includegraphics[width=10cm]{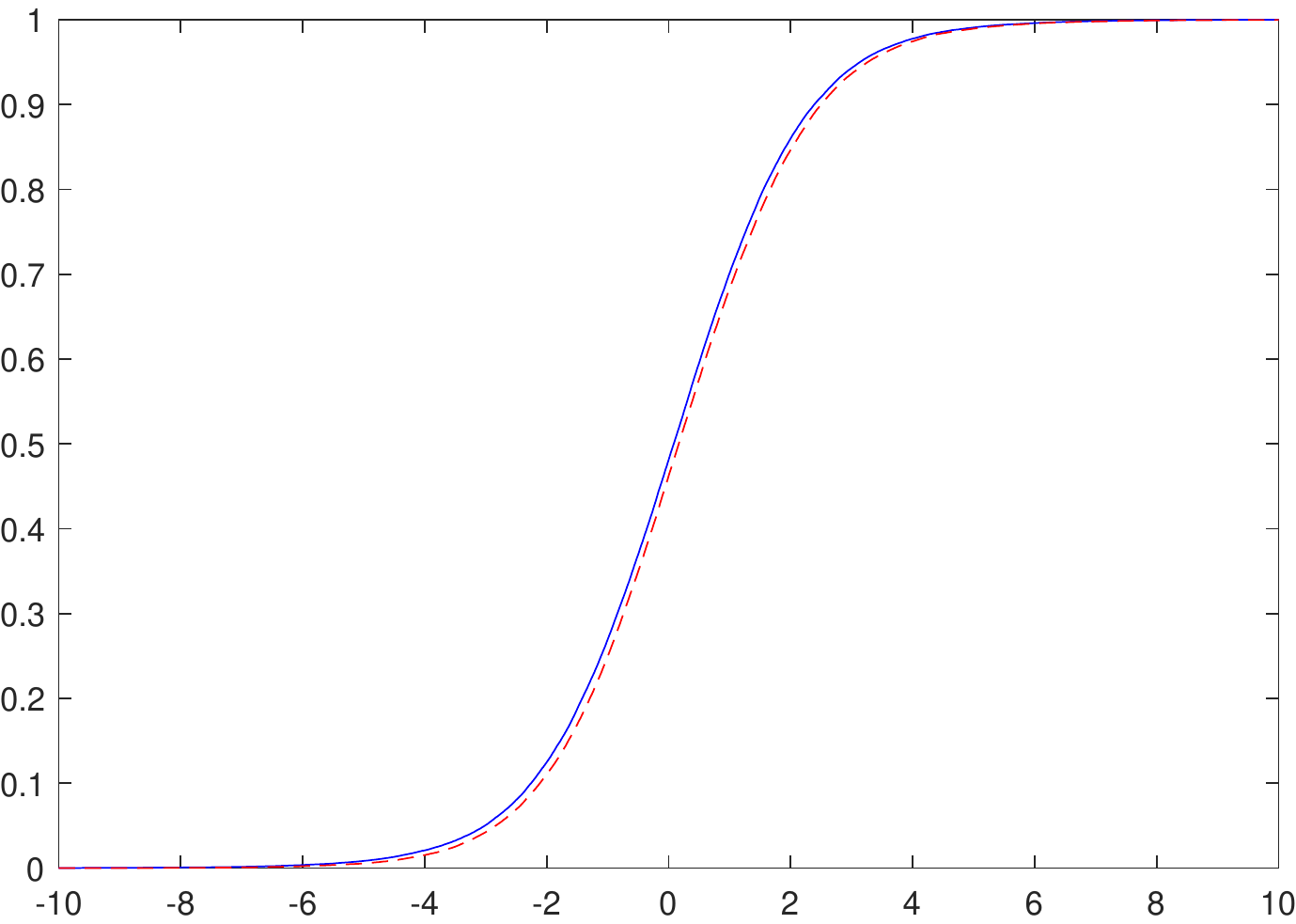}
  \caption{Empirical cdf of $n^{1/2}(\hat\alpha_{n,k_n^*}-\alpha)$ for $n=2,000$ (red, dashed) and limit cdf according to \eqref{eq:HillParbreakasymp} (blue, solid)}
  \label{fig:HillcompParbreak}
\end{figure}

Figure \ref{fig:HillcompParbreak} displays the normalized estimation error $n^{1/2}(\hat\alpha_{n,k_n^*}-\alpha)$ of the Hill estimator for sample size $n=2,000$ and its asymptotic approximation given in \eqref{eq:HillParbreakasymp}. As in the Pareto case, the approximation is more accurate than that for $k_n^*/n$, but it is slightly worse than in the Pareto case. Again the underestimation of the break point leads to a substantial increase of the estimation error compared with the best possible choice of $k$; the RMSE of $\hat\alpha_{n,k_n^*}$ is about 31\% higher than the minimal value which is achieved for $k=967$ and about 23\% higher than that of $\hat\alpha_{n,nt_0}$.

\begin{figure}[tb]
  \centering

  \includegraphics[width=10cm]{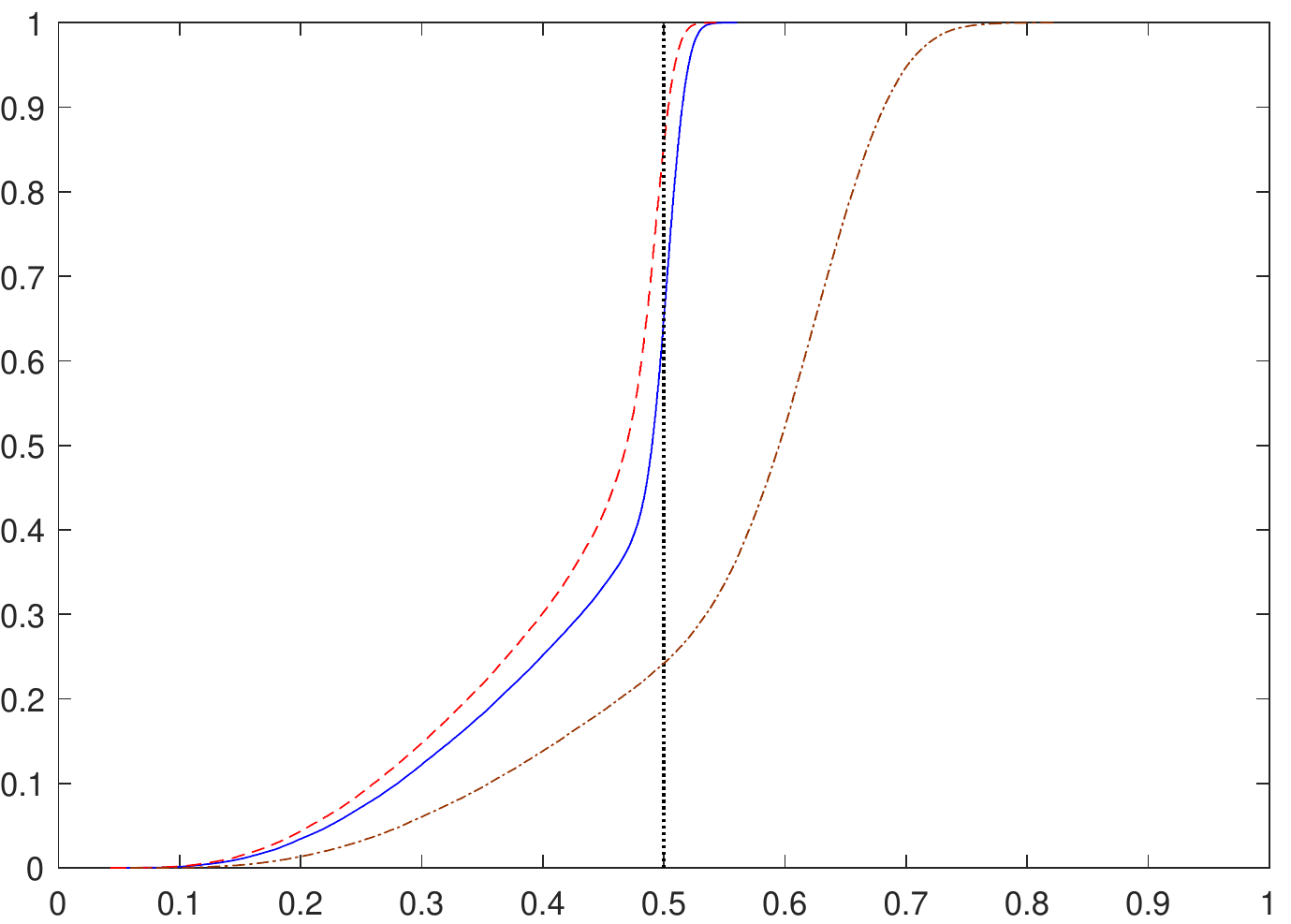}
  \caption{Empirical cdf of $k_n^*/n$ for the models described by $H_c$ (blue, solid), $H_j$ (red, dashed) and $H_d$ (brown, dash-dotted) with sample size $n=2,000$; the break point $t_0=1/2$ is indicated by the vertical line.}
  \label{fig:optkcompbreakseveralH}
\end{figure}

Finally, we briefly discuss simulation results for models where either $H$ is not continuous at $t_0$ or differentiable at $t_0$. More concretely, in addition to the aforementioned function $H_c(t)=H(t)=\big((t/t_0)^{1/\alpha-1/\beta}-1\big)1_{[t_0,1)}(t)$, we consider the functions $H_j(t)=\big((t/t_0)^{1/\alpha-1/\beta}-2\big)1_{[t_0,1)}(t)$ (i.e., the function jumps by $-1$ at $t_0$ and behaves otherwise like $H_c$) and $H_d(t)=-(t-t_0)^21_{[t_0,1)}(t)$. Figure \ref{fig:optkcompbreakseveralH} shows the empirical distribution of $k_n^*/n$ for these three different models. As expected, in the smooth model $H_d$, the MDSP is not able to detect the structural change at $t_0$ quickly. Moreover, due to the specific effect of the deviation from the Pareto model in the neighborhood of $t_0$, the procedure detects the structural break less accurately in the discontinuous model $H_j$ than in the continuous, but not differentiable model given by $H_c$, which puts more mass in the neighborhood of the break point $t_0=1/2$ (indicated by the dotted vertical line). It is worth mentioning that these differences are not reflected in the relative performance of the Hill estimator: for $H_j$ and $H_d$ the increase of the RMSE of $\hat\alpha_{n,k_n^*}$ relative to the Hill  estimator with minimal RMSE equals 30\%, respectively 31\%.

\begin{remark}{\rm
  In extreme value theory, it is often not assumed that above some
  threshold the tail is exactly of Pareto type, but that the
  difference between the actual tail and the approximating Pareto
  vanishes as the threshold increases. More precisely, a so-called
  second order condition may be used, e.g., that
  $F^\leftarrow(1-tx)/F^\leftarrow(1-t)-x^{1/\alpha}\sim t^\rho g(x)$
  as $t\downarrow 0$ for some non-degenerate function $g$ and some $\rho>0$. In such a
  setting, a plethora of methods for selecting $k$ aiming at a minimal
  RMSE of the Hill estimator have been suggested; see, for instance,
  \cite{Beiretal04}, Section 4.7, or \cite{GO01} for a comparison of
  some of these procedures.

  Using the approach employed in the proofs of Theorem \ref{th:Parmain} and Theorem \ref{th:Parbreak}, one can analyze the behavior of the MDSP in such a framework, too, provided the minimum of the Kolmogorov-Smirnov distance $D_k$ is considered only over a set of indices $k\in\{2,\ldots, k_n\}$ for some so-called intermediate sequence $k_n$, i.e., $k_n\to\infty$, but $k_n/n\to 0$. (A related result can be found in \cite{KP08} which considers the Kolmogorov-Smirnov distance for an intermediate sequence converging to $\infty$ sufficiently slowly such that the deviation from the Pareto model is asymptotically negligible.) It turns out that, similarly as in the case of a differentiable function $H$, this procedure is not able to pick a value $k$ that asymptotically minimizes the RMSE of the Hill estimator. Moreover, simulations show that in terms of the RMSE of the Hill estimator, it is usually outperformed by other methods like the sequential procedure proposed by \cite{DK98} or the bootstrap approach examined by \cite{DdHPdV01}.}

\end{remark}

\section{Linear preferential attachment (PA) networks}
  \label{sec:networks}

So far we have only discussed the performance of the MDSP for iid data. However, one important application of the MDSP are network models, where the power-law behavior of the degree distribution is widely observed.
Theoretically, the linear PA model asymptotically generates power-law degree distributions and is therefore a popular choice to
model networks. In this section, we first give an overview of the linear PA model and discuss the tail behavior of the in- and out-degrees, and then summarize simulation results on the performance of
the MDSP for such dependent data.

\subsection{The linear PA model}\label{subsec:linpref}
The directed edge linear PA model
\cite{BBCR03,KR01}
constructs  a growing directed random graph $G(n)=(V(n),E(n))$ whose
dynamics depend on five non-negative real numbers
$\alpha,
\beta, \gamma$, $\delta_{\text in}$ and $\delta_{\text out}$, where
$\alpha+\beta+\gamma=1$ and $\din,\dout
>0$. The values of
  $\alpha, \beta$ and $\gamma$ are probabilities for different
  scenarios of graph growth.
The purpose of  $\delta_{\text in}$ and
  $\delta_{\text out}$  as parameters is to either de-emphasize or
  emphasize the importance of degrees in the specification of
  attachment probabilities if the $\delta$'s are large or small, respectively.
   Setting either
   $\delta_\text{in}$ or $\delta_\text{out}$ to zero makes nodes
   created without in- or out-degrees remain without in- or
  out-degrees; this is unrealistic for many applications.
So here we assume both $\din,\dout >0$ and the detailed model construction follows.

To avoid degenerate situations,
assume  that each of the numbers $\alpha,
\beta, \gamma$ is strictly smaller than 1.
We obtain a new graph $G({n})$
by adding  one edge to the existing graph $G({n-1})$ and
index the constructed graphs by the number $n$ of edges in $E(n)$.
We start with an arbitrary initial finite directed
graph $G({n_0})$ with at least one node and $n_0$ edges.
For $n >n_0$,
$G(n)=(V(n),E(n))$ is a graph with $|E(n)|=n$ edges and a random number $|V(n)|=N(n)$ of
nodes. If $u\in V(n)$, $D_{\rm in}^{(n)}(u)$ and $D_{\rm out}^{(n)}(u)$
denote the in- and out-degree of $u$ respectively in $G(n)$.
There are three
scenarios that we call the $\alpha$, $\beta$ and
$\gamma$-schemes, which are activated by flipping a
3-sided coin whose outcomes are $1,2,3$ with probabilities
$\alpha,\beta,\gamma$. More formally,  we have an iid sequence
 of multinomial random
variables $\{J_n, n>n_0\}$ with cells labelled $1,2,3$ and cell
probabilities $\alpha,\beta,\gamma$.
Then the graph $G(n)$ is
obtained from  $G(n-1)$ as follows.

\tikzset{
    >=stealth',
    punkt/.style={
           rectangle,
           rounded corners,
           draw=black, very thick,
           text width=6.5em,
           minimum height=2em,
           text centered},
    pil/.style={
           ->,
           thick,
           shorten <=2pt,
           shorten >=2pt,}
}
\newsavebox{\mytikzpic}
\begin{lrbox}{\mytikzpic}
     \begin{tikzpicture}
    \begin{scope}[xshift=0cm,yshift=1cm,scale=.8]
      \node[draw,circle,fill=white] (s1) at (2,0) {$v$};
      \node[draw,circle,fill=gray!30!white] (s2) at (.5,-1.5) {$w$};
      \draw[->] (s1.south west)--(s2.north east){};
      \draw[dashed] (0,-2.2) circle [x radius=2cm, y radius=15mm];
    \end{scope}

     \begin{scope}[xshift=4cm,yshift=1cm,scale=.8]
      \node[draw,circle,fill=gray!30!white] (s1) at (.5,-1.5) {$v$};
      \node[draw,circle,fill=gray!30!white] (s2) at (-.5,-2.5) {$w$};
      \draw[->] (s1.south west)--(s2.north east){};
      \draw[dashed] (0,-2.2) circle [x radius=2cm, y radius=15mm];
    \end{scope}

     \begin{scope}[xshift=8cm,yshift=1cm,scale=.8]
      \node[draw,circle,fill=white] (s1) at (2,0) {$w$};
      \node[draw,circle,fill=gray!30!white] (s2) at (.5,-1.5) {$v$};
      \draw[->] (s2.north east)--(s1.south west){};
      \draw[dashed] (0,-2.2) circle [x radius=2cm, y radius=15mm];
    \end{scope}

     \node at (0,-3.5) {$\alpha$-scheme};
     \node at (4,-3.5) {$\beta$-scheme};
     \node at (8,-3.5) {$\gamma$-scheme};
  \end{tikzpicture}

\end{lrbox}
  \begin{figure}[h]
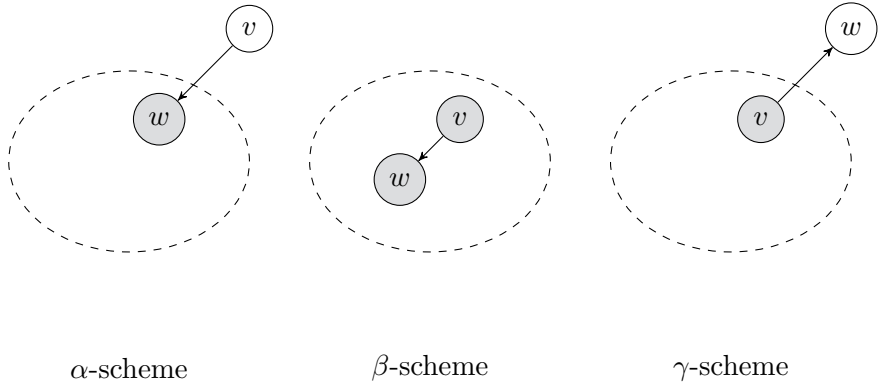

    \centering
    \usebox{\mytikzpic}
    \caption{Evolution of a linear preferential attachment model.}
\end{figure}

\begin{itemize}
\item
If $J_n=1$ (with probability
$\alpha$),  append to $G(n-1)$ a new node $v\in V(n)\setminus V(n-1)$ and an edge
$(v,w)$ leading
from $v$ to an existing node $w \in V(n-1)$.
Choose the existing node $w\in V(n-1)$  with probability depending
on its in-degree in $G(n-1)$:
\beqq \label{eq:probIn}
\PP[\text{choose $w\in V(n-1)$}] \,=\, \frac{D_{\rm in}^{(n-1)}(w)+\delta_{\text
    in}}{n-1+\delta_{\text in}N(n-1)} \,.
\eeqq
\item If $J_n=2$ (with probability $\beta$), add a directed edge
$(v,w) $ to $E({n-1})$ with $v\in V(n-1)=V(n) $ and $w\in V(n-1)=V(n) $ and
 the existing nodes $v,w$ are chosen independently from the nodes of $G(n-1)$ with
 probabilities
$$ \label{eq:probBoth}
\PP[\text{choose $(v,w)$}] \,=\, \Bigl(\frac{D_{\rm out}^{(n-1)}(v)+\delta_{\text
    out}}{n-1+\delta_{\text out}N(n-1)}\Bigr)\Bigl(
 \frac{D_{\rm in}^{(n-1)}(w)+\delta_{\text
    in}}{n-1+\delta_{\text in}N(n-1)}\Bigr).
$$
\item  If $J_n=3$ (with probability
$\gamma$),  append to $G(n-1)$ a new node $w\in V(n)\setminus V(n-1)$ and an edge $(v,w)$ leading
from  the existing node $v\in V(n-1)$  to the new node $w$. Choose
the existing node $v\in V(n-1)$ with probability
\beqq \label{eq:probOut}
\PP[\text{choose $v \in V(n-1)$}] \,=\, \frac{D_{\rm out}^{(n-1)}(v)+\delta_{\text
    out}}{n-1+\delta_{\text out}N(n-1)}\,.
\eeqq
\end{itemize}
{
Note that this construction allows the possibility of having self loops in the case where $J_n=2$,
but the proportion of edges that are self loops goes to 0 as $n\to\infty$. Also, multiple edges are allowed between two nodes.}

\subsection{Power law of degree distributions}

Given an observed network with $n$ edges, let $N_{ij}(n)$ denote the
number of nodes in $G(n)$ with in-degree $i$ and out-degree $j$. If the network
is generated from the linear PA model described above, then from
\cite{BBCR03},
  there exists a proper probability distribution $\{f_{ij}\}$ such
  that almost surely
\beqq\label{pij}
\frac{N_{ij}(n) }{N(n)} \,\to\, f_{ij},\quad n\to\infty.
\eeqq
Thus, the empirical frequency of nodes with in- and out-degrees $(i,j)$ converges to a limiting, non-random, theoretical frequency.
Consider the limiting marginal in-degree distribution $\fin_i:=\sum_{j=0}^\infty f_{ij}$.
From \cite[Theorem 3.1]{BBCR03},
\beqq
\fin_i \sim \Cin i^{-(1+\ain)}\mbox{ as }i\to\infty, \quad\text{as long as } \alpha\din+\gamma>0,\label{asyI}
\eeqq
for some finite positive constant $\Cin$,
and the power-law index
\beqq\label{c1}
\ain \,=\, \frac{1+\din(\alpha+\gamma)}{\alpha+\beta}.
\eeqq
Similarly, the limiting marginal out-degree distribution has the same property:
\begin{align*}
\fout_j &\,:=\, \sum_{i=0}^\infty f_{ij} \sim \Cout j^{-(1+\aout)}\mbox{ as }j\to\infty, \quad\text{as long as }\gamma\dout+\alpha>0,
\end{align*}
for some $\Cout$  positive and
\beqq\label{c2}
\aout \,=\, \frac{1+\dout(\alpha+\gamma)}{\beta+\gamma}.
\eeqq
Limit behavior of the degree counts in this linear PA model is studied in
\cite{krapivsky:2001,  resnick:samorodnitsky:2016b, resnick:samorodnitsky:towsley:davis:willis:wan:2016, wang:resnick:2016, wang:resnick:2015}.

Two parametric estimation methods for this directed linear PA model are derived in \cite{wan:wang:davis:resnick:2017}, giving
estimates of $\ain$ and $\aout$ by simply plugging in the estimated parameters into \eqref{c1} and \eqref{c2}, respectively.
However, these estimates rely heavily on the correctness of the underlying model, which is hard to guarantee for real data.
In \cite{wan:wang:davis:resnick:2019}, another estimation method coupling the Hill estimation of marginal degree
distribution tail indices with the MDSP  is proposed.
Despite the dependence structure of degree
  sequences, the consistency of the Hill estimator is proved for
  certain linear PA models in
  \cite{wang:resnick:2018a,wang:resnick:2018b}, provided that the
  sequence $k_n$ of order statistics used for estimation is
  deterministic and increases at a suitable rate. In
  \cite{wan:wang:davis:resnick:2019}, the performance of parametric
  and Hill
  approaches to estimating the tail indices is compared via
  simulation.
  The MLE parametric approach  is more efficient in the absence of
  model error but the Hill estimator combined with the MDSP
  is much more robust against deviations from the linear PA model or
  data corruption. Replacing the MDSP with some better selection
  procedure may hopefully increase the efficiency of the Hill
  estimator without compromising its robustness.

\subsection{Simulations}
\label{subsec:sim}
We now further examine the performance of the MDSP in the context of linear PA models through simulations.
We simulate 10,000 linear PA graphs with expected number $m=10^6$ of edges.
To this end,
starting from a trivial core with just one node and no edges, we grow the network as described above until it has $n=\ceil{(\alpha+\gamma)m}$ nodes. Then we try to find an appropriate number $k$ such that the distribution of the $k$ largest observed in-degrees can be well fitted by a power tail.

To work with realistic models, we have chosen the parameters of the simulated linear PA network equal to the estimates obtained from real networks using the snap shot methodology
described in \cite{wan:wang:davis:resnick:2017}. Here we report the results for the following two parameter sets estimated from
KONECT \cite{kunegis:2013}
data sets:
\begin{description}
\item[Example I:] 
  For {\em Baidu related pages}
  (\url{http://konect.uni-koblenz.de/networks/zhishi-baidu-relatedpages
  }) we obtained the estimates $(\alpha,\beta,\gamma,\din,\dout)=(0.0978,0.873,0.0289,2.05,$ $13.8)$
   resulting in $\ain=1.30$;
\item[Example II:]  {\em
  Facebook wall posts}
(\url{http://konect.uni-koblenz.de/networks/facebook-wosn-wall}) \linebreak
leads to the estimates  $(\alpha,\beta,\gamma,\din,\dout)=(0.0327,0.946,0.0209,8.88,9.59)$
   which implies $\ain=1.51$.
\end{description}
(The parameters are rounded to 3 significant digits.)
\smallbreak

\begin{figure}[tb]
  \centering

  \includegraphics[width=13cm]{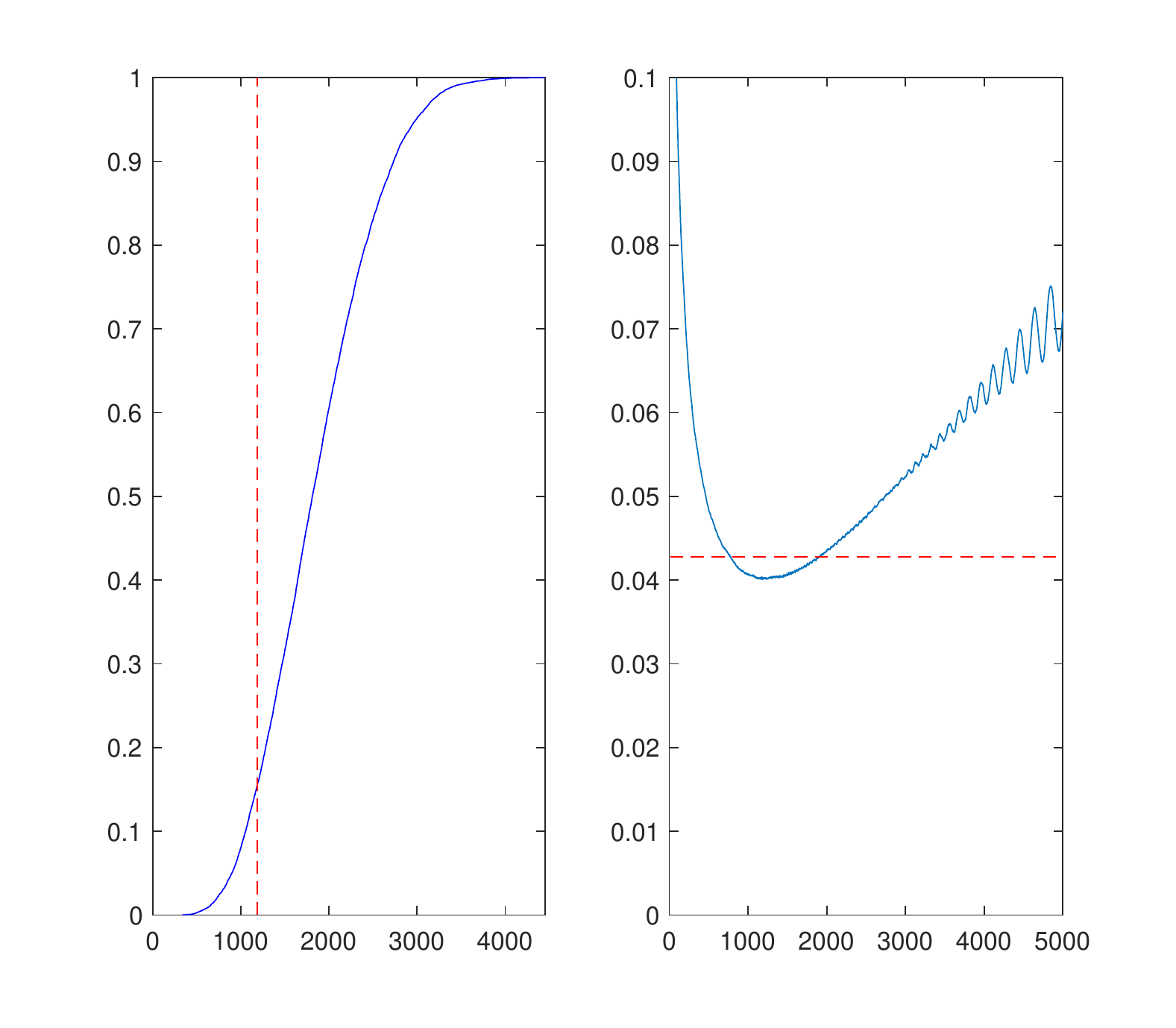}
  \caption{Left: empirical cdf of $k_n^*$ for the linear PA Model I, the RMSE minimizing value of $k$ is indicated by the dashed red line; right: RMSE of the Hill estimator vs.\ $k$, the RMSE of $\hat\alpha_{n,k_n^*}$ is indicated by the dashed red line.}
  \label{fig:LPAN1}
\end{figure}

In Example I the RMSE of the Hill estimator $\hat\alpha_{n,k_n^*}$ for the tail index of the in-degree is just 6.8\% larger than the minimal RMSE over all deterministic choices of $k\in\{10,\ldots,10000\}$ which is attained for $k=1187$. In this respect, the MDSP works much better for this linear PA network model than in any situation with iid data considered in Section \ref{sec:Pareto} and Section \ref{sec:Parbreak}. The right plot of Figure \ref{fig:LPAN1}, which shows the RMSE of the Hill estimator as a function of $k$, hints at the reason for this good performance. There is a wide range of values $k$ that lead to almost the same RMSE. So although the distribution of $k_n^*$ (shown in the left plot) is spread out over the interval $[500,3500]$, this does not increase the RMSE substantially.

\begin{figure}[tb]
  \centering

  \includegraphics[width=13cm]{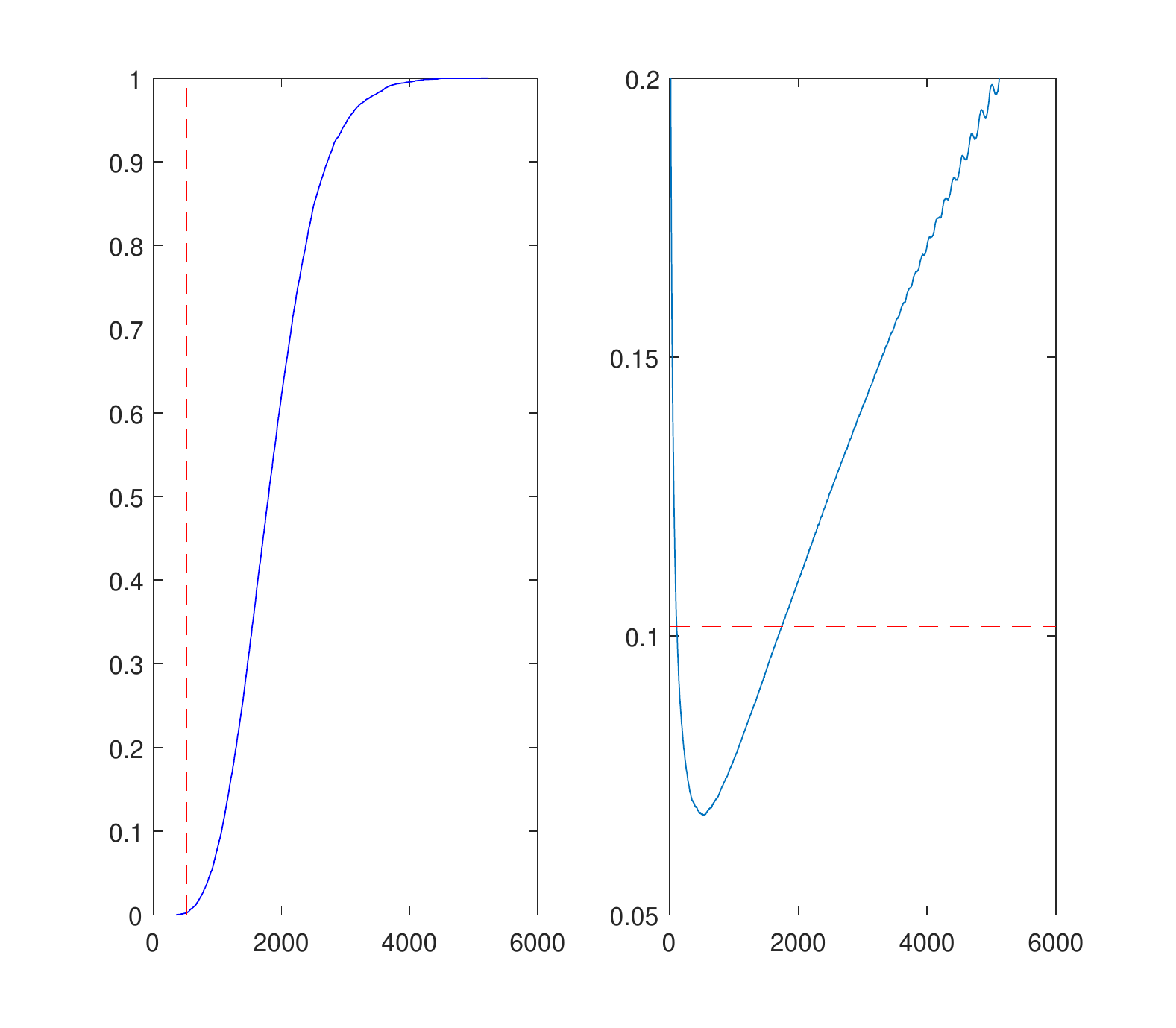}
  \caption{Left: empirical cdf of $k_n^*$ for the linear PA Model II, the RMSE minimizing value of $k$ is indicated by the dashed red line; right: RMSE of the Hill estimator vs.\ $k$, the RMSE of $\hat\alpha_{n,k_n^*}$ is indicated by the dashed red line.}
  \label{fig:LPAN2}
\end{figure}

In the
Facebook Example II, the loss of efficiency is much larger. Here the RMSE of $\hat\alpha_{n,k_n^*}$ is about $50.0\%$ larger than the minimal RMSE. For this model, the RMSE increases much faster as $k$ deviates from the RMSE-minimizing value  $k=523$ (see the right plot of Figure \ref{fig:LPAN2}). Since the distribution of $k_n^*$ (with an estimated mean of 1857) puts almost all its mass on values of $k$ much larger than the point of minimum, the sensitivity of the Hill estimator to an inappropriate selection of $k$ leads to a rather poor performance of  $\alpha_{n,k_n^*}$.

Though the performance of the MDSP in our simulation is somewhat mixed, it yields good results in terms of the RMSE of the Hill estimator if the Hill estimator is not very sensitive to the choice of the threshold. According to further simulation results (not reported here), such a behavior seems to be more common for network data than for many popular models of iid data.
Hence, we conclude that the MDSP often works well on the linear PA models under proper choices of parameters.

\section{Conclusions}

We discussed the asymptotic and the finite sample performance of the minimum distance selection procedure. It was shown for models of iid data with Pareto tail that, unlike previously proposed methods, the sample fraction $k_n^*/n$ chosen by the MDSP does not asymptotically concentrate on one point. Instead, it often yields too small a value of $k$ if there is a clear structural break in the distribution below some threshold, leading to a strongly increased variance and RMSE of the Hill estimator. On the other hand, if there is a smooth transition from the Pareto tail to a different behavior for smaller observations, the MDSP usually picks up this change point too late, again leading to a substantial increase of the RMSE of the Hill estimator.

As the simulations of dependent data from linear preferential attachment networks have
shown, the spread of the distribution of $k_n^*$ need not always
result in a large loss of efficiency of the Hill estimator. This is
particularly true if the Hill estimator is rather insensitive to the
choice of $k$ over a wide range of $k$, because the increase in the bias with growing $k$ is balanced by the decrease of the variance.

There may be another reason why the MDSP shows a considerably different behavior for the in-degrees of linear preferential attachment networks and for iid observations. In the latter situation, for any fixed sample size, all observations are drawn from a distribution which is regularly varying, that is $(1-F(tx))/(1-F(x))\to t^{-1/\alpha}$ as $x\to\infty$. This is not true for the in-degrees of a linear PA model with a given deterministic number of nodes, which cannot take on values larger than the number of edges in the network. In fact, the distribution of the in-degrees changes when the sample size (i.e., the number of nodes) increases, and only in the limit (as described by \eqref{pij} and \eqref{asyI}) the distribution has a power tail behavior. So while in the situation considered in the Sections \ref{sec:Pareto} and \ref{sec:Parbreak} the tail index $\alpha$ has the same operational meaning for each sample size, for linear PA models it is defined only as a limit parameter as the number $n$ of nodes tends to infinity. Since, strictly speaking, for any fixed $n$, there is no tail index, the interpretation of the RMSE of the Hill estimator is somewhat unclear in this setting.

\appendix
\section{Proofs}

\begin{proofof} Theorem \ref{th:Parmain}.\rm\quad
Let $U_i$, $i\in\N$, be iid uniform rv's, $\xi_i$, $i\in\N$, be iid standard exponential rv's, and $S_k:=\sum_{i=1}^k \xi_i$. Then $(U_{k:n})_{1\le k\le n}=^d (S_k/S_{n+1})_{1\le k\le n}$ for all $n\in\N$ (\cite{Reiss89}, Cor.\ 1.6.9). Hence, by the quantile transformation, it suffices to prove that the assertion holds for
\begin{equation} \label{eq:Dktdef}
 \tilde D_k = \max_{1\le j\le k} \bigg|\Big(\frac{F^\leftarrow(1-S_j/S_{n+1})}{F^\leftarrow(1-S_k/S_{n+1})}\Big)^{-\tilde\alpha_k}-\frac jk\bigg|
\end{equation}
instead of $D_k$ with
$$\tilde\alpha_{n,k} := \bigg(\frac 1{k-1}\sum_{i=1}^{k-1} \log\frac{F^\leftarrow(1-S_i/S_{n+1})}{F^\leftarrow(1-S_k/S_{n+1})}\bigg)^{-1}
$$
and $F^\leftarrow$ denoting the quantile function pertaining to $F$.
For $F$ according to \eqref{eq:Parcdf} this simplifies to
\begin{equation} \label{eq:DktPar}
 \tilde D_k = \max_{1\le j\le k} \bigg|\Big(\frac{S_j}{S_k}\Big)^{\tilde\alpha_{n,k}/\alpha}-\frac jk\bigg|
\end{equation}
with
\begin{equation} \label{eq:alphatildefrac}
 \frac \alpha{\tilde\alpha_{n,k}} = \frac 1{k-1}\sum_{i=1}^{k-1} \log \frac{S_k}{S_i} = \frac 1{k-1}\sum_{i=1}^{k-1} \log \frac ki - \frac 1{k-1}\sum_{i=1}^{k-1} \log \frac{kS_i}{iS_k}.
\end{equation}
Note that neither $\tilde D_k$ nor $\tilde\alpha_{n,k}$ depend on $n$, so that we will drop the index $n$ when using the latter in the remaining part of the proof.

The first sum on the right hand side of \eqref{eq:alphatildefrac} is a Riemann approximation of $\int_0^1\log(1/t)\, dt$ $=1$ with an approximation  error $O((\log k)/k)$. To analyze the second sum, we use the so-called Hungarian construction (see \cite{KMT75,KMT76}): for suitable versions of the $\xi_i$, there exists a Brownian motion $W$ such that
$$ \max_{1\le i\le k} |S_i-i-W(i)| = O(\log k)\quad \text{a.s.} $$
Let $LL_x:= \log\log(\e^\e\vee x)$. Then
\begin{align} \label{eq:ratioapprox}
  \frac{kS_i}{iS_k}-1 & = \frac{W(i)-(i/k)W(k)+O(\log k)}{i+(i/k)W(k)+O((i/k)\log k)} \nonumber\\
   & = \frac{W(i)}i-\frac{W(k)}k + O\Big(\frac{\log k}i + \Big(\frac{LL_i LL_k}{ik}\Big)^{1/2}\Big)\nonumber\\
   & = \frac{W(i)}i-\frac{W(k)}k + O\Big(\frac{\log k}i\Big)
\end{align}
uniformly for all $1\le i\le k$, where in the second step we have used the law of iterated logarithm.

It follows by the strong law of large numbers and a Taylor expansion of $\log$ that
\begin{align*}
\frac 1{k-1}\sum_{i=1}^{k-1} \log \frac{kS_i}{iS_k}
 & =   \frac 1{k-1}\sum_{i=\ceil{\log k}}^{k-1} \log\bigg[ 1+\frac{W(i)}i-\frac{W(k)}k + O\Big(\frac{\log k}i\Big)\bigg] + O\Big(\frac{\log k}k\Big)\\
 & =  \frac 1{k-1}\sum_{i=\ceil{\log k}}^{k-1} \Big(\frac{W(i)}i-\frac{W(k)}k +O\Big(\frac{\log k}i\Big)\Big)+ O\Big(\frac{\log k}k\Big)\\
 & =  \frac 1{k-1}\sum_{i=\ceil{\log k}}^{k-1} \Big(\frac{W(i)}i-\frac{W(k)}k\Big)+ O\Big(\frac{\log^2 k}k\Big),
\end{align*}
since $\sum_{i=1}^k i^{-1}=O(\log k)$.
Moreover, by the law of iterated logarithm,
$$  \frac 1{k-1}\sum_{i=1}^{\ceil{\log k}-1}\Big( \frac{W(i)}i-\frac{W(k)}k\Big) 
= O\Big(\frac{\log k}k\Big).
$$
To sum up, we have shown that
$$ \frac \alpha{\tilde\alpha_k}-1 = -\frac 1{k-1}\sum_{i=1}^{k-1}\Big( \frac{W(i)}i-\frac{W(k)}k\Big)+O\Big(\frac{\log^2 k}k\Big),
$$
which in turn implies
\begin{equation} \label{eq:alphatildeapprox}
  \frac{\tilde\alpha_k}\alpha =1 +\frac 1{k-1}\sum_{i=1}^{k-1} \Big(\frac{W(i)}i-\frac{W(k)}k\Big)+O\Big(\frac{\log^2 k}k\Big).
\end{equation}

Let $\tau_k=\tilde\alpha_k/\alpha-1 = O((LL_k/k)^{1/2})$.
Then, by \eqref{eq:ratioapprox}, one has, uniformly for all $1\le j\le k$,
$$ \Big(\frac{S_j}{S_k}\Big)^{\tilde\alpha_k/\alpha} = \Big(\frac jk\Big)^{1+\tau_k} \Big[1+\frac{W(j)}j-\frac{W(k)}k+O\Big(\frac{\log k}j\Big)\Big]^{1+\tau_k}. $$
The first factor on the right hand side equals
$$ \frac jk\Big( 1+\tau_k\log\frac jk + O\Big(\frac{\log^2(j/k)LL_k}k\Big)\Big). $$
A Taylor expansion of $\log$ and $\exp$ shows that the second factor is equal to
\begin{eqnarray*}
  \lefteqn{\exp\Big((1+\tau_k)\log\Big[1+\frac{W(j)}j-\frac{W(k)}k+O\Big(\frac{\log k}j\Big)\Big]\Big)}\\
  & = & \exp\Big((1+\tau_k)\Big[\frac{W(j)}j-\frac{W(k)}k+O\Big(\frac{\log k}j\Big)\Big]\Big)\\
  & = & 1+\frac{W(j)}j-\frac{W(k)}k+O\Big(\frac{\log k}j\Big),
\end{eqnarray*}
because $\tau_k (LL_j/j)^{1/2} = o((\log k)/j)$.
Therefore, since $t(\log t)^2$ is bounded on the unit interval and thus $O(\log^2(j/k)LL_k/k)=o((\log k)/j)$,
\begin{equation} \label{eq:approx1}
   \Big(\frac{S_j}{S_k}\Big)^{\tilde\alpha_k/\alpha} = \frac jk \Big(1+ \frac{W(j)}j-\frac{W(k)}k+\tau_k\log\frac jk + O\Big(\frac{\log k}j\Big)\Big)
\end{equation}
uniformly for all $1\le j\le k$.
Combining \eqref{eq:DktPar}, \eqref{eq:alphatildeapprox} and \eqref{eq:approx1}, we arrive at
\begin{equation} \label{eq:Dktapprox}
  \tilde D_k  = \max_{1\le j\le k} \bigg|\frac jk \Big(\frac{W(j)}j-\frac{W(k)}k\Big) + \frac jk \log\frac jk \Big(\frac 1{k-1}\sum_{i=1}^{k-1} \frac{W(i)}i-\frac{W(k)}k\Big)\bigg|+O\Big(\frac{\log^2 k}k\Big).
\end{equation}

In the last step, we replace the maximum over the discrete points $j$ with a supremum over a whole interval and the sum with an integral. To this end, for each $n$, we define the Brownian motion $W_n(x)=n^{-1/2}W(nx)$, $x\ge 0$. Then, with $k=\ceil{nt}$ and $j=sk$
\begin{align} \label{eq:Wnreplace}
  n^{1/2} \tilde D_{\ceil{nt}} & = \max_{\substack{s\in (0,1]\\s\ceil{nt}\in\N}} \bigg| s\Big( \frac{W_n(s\ceil{nt}/n)}{s\ceil{nt}/n} - \frac{W_n(\ceil{nt}/n)}{\ceil{nt}/n}\Big)+ \nonumber\\
  & \hspace{1cm}
  + s\log s\Big( \frac 1{\ceil{nt}-1}\sum_{i=1}^{\ceil{nt}-1} \frac{W_n(i/n)}{i/n} - \frac{W_n(\ceil{nt}/n)}{\ceil{nt}/n}\Big)\bigg|
  + O\Big(\frac{\log^2(nt)}{n^{1/2}t}\Big).
\end{align}
Recall that the modulus of continuity of a Brownian motion on the unit interval equals $\omega_W(\delta)=(2\delta|\log\delta|)^{1/2}$ a.s. Hence
\begin{equation} \label{eq:Wncont1}
  \sup_{s\in (0,1]} \big| W_n(s\ceil{nt}/n)-W_n(st)\big| = O_P\Big(\Big(\frac{\log n}n\Big)^{1/2}\Big).
\end{equation}
Furthermore, for all $1/n\le x\le y\le 1$
\begin{equation} \label{eq:Wncont2}
  \Big|\frac{W_n(y)}y-\frac{W_n(x)}x\Big| \le \frac{|W_n(y)-W_n(x)|}y+|W_n(x)|\frac{|y-x|}{xy}.
\end{equation}
Thus, we conclude using the law of iterated logarithm (at 0) that, uniformly for $t\in[2/n,1], s\in[1/\ceil{nt},1]$,
\begin{align} \label{eq:Wncont3}
 s\Big| \frac{W_n(s\ceil{nt}/n)}{s\ceil{nt}/n}-\frac{W_n(st)}{st}\Big|
 & = O_P\Big(\Big(\frac{\log n}n\Big)^{1/2} \frac 1t + (stLL_{1/(st)})^{1/2}\frac{s^2/n}{(st)^2}\Big)\nonumber\\
 & = O_P\Big(\Big(\frac{\log n}n\Big)^{1/2}\frac 1t\Big).
\end{align}
Finally, again using \eqref{eq:Wncont2}, we obtain
 \begin{align} \label{eq:intapprox}
   \frac 1{\ceil{nt}-1} & \sum_{i=1}^{\ceil{nt}-1} \frac{W_n(i/n)}{i/n} - \int_0^1 \frac{W_n(tx)}{tx}\, dx \nonumber\\
   & =  \sum_{i=1}^{\ceil{nt}-1} \int_{(i-1)/\ceil{nt-1}}^{i/\ceil{nt-1}} \frac{W_n(i/n)}{i/n}- \frac{W_n(tx)}{tx}\, dx \nonumber\\
   & =  O_P\bigg(\frac 1{\ceil{nt}-1}\bigg[ \sum_{i=1}^{\ceil{nt}-1} \frac ni \Big(\Big|\log\frac{t}{\ceil{nt}-1}\Big|\frac{t}{\ceil{nt}-1}\Big)^{1/2} + \nonumber\\
   & \hspace*{4cm} +\Big(LL_{n/i}\frac in\Big)^{1/2}\frac ni \frac t{\ceil{nt}-1}\bigg]\bigg) \nonumber\\
  & =  O_P\Big(\Big(\frac{\log n}n\Big)^{1/2}  \frac{\log(nt)}t \Big)
 \end{align}
 uniformly for $t\in[2/n,1]$.
 Combining \eqref{eq:Wnreplace}, \eqref{eq:Wncont1}, \eqref{eq:Wncont3} and \eqref{eq:intapprox}, we conclude
 \begin{align} \label{eq:Dktapprox2}
   n^{1/2} \tilde D_{\ceil{nt}} & = \max_{\substack{s\in (0,1]\\s\ceil{nt}\in\N}} \bigg| s\Big( \frac{W_n(st)}{st} - \frac{W_n(t)}{t}\Big)+ s\log s\Big( \int_0^1 \frac{W_n(tx)}{tx}\, dx - \frac{W_n(t)}t\Big)\bigg| \nonumber\\
  & \hspace{3cm}
  +  O_P\Big(\frac{\log(nt)(\log(nt)+(\log n)^{1/2})}{n^{1/2}t}\Big)
\end{align}
uniformly for $t\in[2/n,1]$.

To replace the maximum with a supremum over all $s\in(0,1]$, observe that for any $s\in(0,1]$ there is a point $\tilde s$ with $|s-\tilde s|<1/(nt)$ that is considered in maximum. Hence by the modulus of continuity of $W_n$, the law of iterated logarithm and the inequality
$|s\log s-\tilde s\log\tilde s|\le |s-\tilde s|(1+|\log(s\wedge \tilde s)|)$, which holds for all $s,\tilde s\in(0,1]$, one has
\begin{align*} 
  \Bigg|& \max_{\substack{s\in (0,1]\\s\ceil{nt}\in\N}} \bigg| s\Big( \frac{W_n(st)}{st} - \frac{W_n(t)}{t}\Big)+ s\log s\Big( \int_0^1 \frac{W_n(tx)}{tx}\, dx - \frac{W_n(t)}t\Big)\bigg| \nonumber\\
  &- \sup_{s\in (0,1]}  \bigg| s\Big( \frac{W_n(st)}{st} - \frac{W_n(t)}{t}\Big)+ s\log s\Big( \int_0^1 \frac{W_n(tx)}{tx}\, dx - \frac{W_n(t)}t\Big)\bigg|\Bigg| \nonumber\\
  & =  O_P\Big(\frac 1t\Big(\frac{\log n}{n}\Big)^{1/2}+\frac 1{nt} \Big(\frac{LL_{1/t}}t\Big)^{1/2}\Big) + O_P\Big(\frac 1{nt}(1+\log(nt))\Big(\frac{LL_{1/t}}t\Big)^{1/2}\Big).
\end{align*}
Now the assertion follows readily.
\end{proofof}

\begin{proofof} Corollary \ref{cor:Par}. \rm
 By Theorem \ref{th:Parmain}, $k^{1/2} D_k\to \sup_{s\in (0,1]} |Z_1(s,1)|$ weakly as $k\to \infty$. According to Skohorod's theorem, there exist versions such that  the convergence holds almost surely. Moreover, all $D_k$ and the limit random variable as well are almost surely strictly positive. Hence, for all sequences $k_n=o(n)$, it follows that $n^{1/2} \min_{2\le k\le k_n} D_k \ge (n/k_n)^{1/2} \min_{2\le k\le k_n} k^{1/2} D_k \to \infty$ almost surely. This implies that with probability tending to 1, $k_n^*$ must be larger than $k_n$.

 Because $k_n=o(n)$ is arbitrary, we conclude that the sequence $n/k_n^*$ is stochastically bounded. Since, by Theorem \ref{th:Parmain},
 $(n^{1/2} D_\ceil{nt})_{t\in[\eps,1]}$ converges weakly to  
 $(Z_1(t))_{t\in[\eps,1]}$ (w.r.t.\ the supremum norm) for all $\eps>0$ and $Z_1$ is continuous on $(0,1]$ with a unique point of minimum, the asymptotic behavior of $k_n^*/n$ follows (cf.\ Corollary 5.58 of \cite{vdV98}). The asymptotics of the Hill estimator can be easily derived from \eqref{eq:alphatildeapprox} and the approximations established in the last part of the proof of Theorem \ref{th:Parmain}, in particular \eqref{eq:Wncont1} and \eqref{eq:intapprox}.
\end{proofof}

\begin{proofof} Theorem \ref{th:Parbreak}. \rm
 Since $n^{-1/2}(\ceil{n(t_0-\eps_n)}-nt_0)\to -\infty$, the central limit theorem yields
 \begin{align*}
  P\{ &X_{n-\ceil{nt}+1:n}  >F^\leftarrow(1-t_0)\}\\
    & = P\Big\{ n^{-1/2}\Big( \sum_{i=1}^{n} 1_{(F^\leftarrow(1-t_0),\infty)}(X_i) -nt_0\Big) \ge n^{-1/2} (\ceil{nt}-nt_0)\Big\} \to 1
 \end{align*}
 uniformly for all $t\in[2/n, t_0-\eps_n]$. Hence assertion (i) is an immediate consequence of Theorem \ref{th:Parmain}.

 To prove the remaining assertions, we use the same approach as in the proof of Theorem \ref{th:Parmain}. In the present setting
 \begin{align*}
   \frac\alpha{\tilde\alpha_{n,k}} & = \frac\alpha{k-1} \sum_{i=1}^{k-1} \log \frac{F^\leftarrow(1-S_i/S_{n+1})}{F^\leftarrow(1-S_k/S_{n+1})} \\
   & = \frac 1{k-1} \sum_{i=1}^{k-1} \log \frac{S_k}{S_i} + \frac\alpha{k-1} \sum_{i=1}^{k-1} \log \frac{1+H(S_i/S_{n+1})}{1+H(S_k/S_{n+1})} \\
   & =: \frac\alpha{\tilde\alpha_k^P} + \Delta_\alpha(n,k).
 \end{align*}
 The first term arises in the Pareto model and has been analyzed in the proof of Theorem \ref{th:Parmain}; in particular, $\alpha/\tilde\alpha_k^P\to 1$ in probability uniformly for all $k\in\{\ceil{n\eps},\ldots,n\}$. For the second term the continuity of $H$ implies, for all $\eps>0$
$$ \sup_{t\in[\eps,1]} \big|\Delta_\alpha(n,\ceil{nt})-\alpha IH(t)\big| \to 0 \quad \text{with}\quad IH(t) := \int_0^t \log\frac{1+H(s)}{1+H(t)}\, ds.
$$
Hence
\begin{equation} \label{eq:alphatildeconv}
   \frac\alpha{\tilde\alpha_{n,\ceil{nt}}} \to 1+\alpha IH(t)
\end{equation}
uniformly for $t\in[\eps,1]$.

In the present setting, the approximative Kolmogorov-Smirnov distance $\tilde D_k$ defined in \eqref{eq:Dktdef} equals
\begin{align} \label{eq:Dktrep}
 \tilde D_k & = \max_{1\le j\le k} \bigg| \Big(\frac{S_j}{S_k}\Big)^{\tilde\alpha_{n,k}/\alpha}\Big(\frac{1+H(S_j/S_{n+1})}{1+H(S_k/S_{n+1})}\Big)^{-\tilde\alpha_{n,k}}-\frac jk\bigg| \nonumber \\
 & = \max_{1\le j\le k} \bigg| \Big(\frac{S_j}{S_k}\Big)^{\tilde\alpha^P_{k}/\alpha}\Big(\frac{S_j}{S_k}\Big)^{(\tilde\alpha_{n,k}-\tilde\alpha^P_{k})/\alpha}\Big(\frac{1+H(S_j/S_{n+1})}{1+H(S_k/S_{n+1})}\Big)^{-\tilde\alpha_{n,k}}-\frac jk\bigg|.
\end{align}
It has been shown in the proof of Theorem \ref{th:Parmain} that for all $\eps>0$
$$ \max_{1\le j\le k}\bigg| \Big(\frac{S_j}{S_k}\Big)^{\tilde\alpha^P_{k}/\alpha}-\frac jk\bigg|=O_P(n^{-1/2}) $$
uniformly for $k\in\{\ceil{n\eps},\ldots,n\}$.
Using \eqref{eq:alphatildeconv}, the strong law of large numbers and the continuity of $H$, we obtain
$$ \sup_{t\in[\eps,1]}\bigg|\tilde D_{\ceil{nt}}-\sup_{s\in(0,1]} s \Big|s^{1/(1+\alpha IH(t))-1}\Big( \frac{1+H(st)}{1+H(t)}\Big)^{-\alpha/(1+\alpha IH(t))}-1\Big|\bigg| \to 0.
$$
The supremum over $s\in(0,1]$ vanishes if and only if
\begin{equation} \label{eq:constcond}
   s^{IH(t)}(1+H(st))/(1+H(t))=1\quad \text{for all } s\in(0,1].
\end{equation}
Because $H(st)=0$  for $s<t_0/t, t>t_0$, \eqref{eq:constcond} implies $H(t)=IH(t)=0$. However, then equation \eqref{eq:constcond} cannot hold for $s$ in a right neighborhood of $t_0/t$ on which $H(st)$ does not vanish.  Since the supremum is a continuous function in $t$, its infimum over $[t_1,1]$ is strictly positive for all $t_1>t_0$. Therefore,
\begin{align} \label{eq:largeDktapprox}
 &n^{1/2}\inf_{t\in[t_1,1]} \tilde D_{\ceil{nt}} \nonumber \\
  & = n^{1/2} \inf_{t\in[t_1,1]}\Big(\sup_{s\in(0,1]} s \Big|s^{1/(1+\alpha IH(t))-1}\Big( \frac{1+H(st)}{1+H(t)}\Big)^{-\alpha/(1-\alpha IH(t))}-1\Big|+o(1)\Big) \to\infty.
\end{align}
By standard arguments one may even find a sequence $t_1=t_{n,1}$ converging to $t_0$ from above sufficiently slowly such that \eqref{eq:largeDktapprox} still holds.

Next we examine $\tilde D_k$ for $k$ such that $k/n$ tends to $t_0$ from above. To this end, we need a refined analysis of $\Delta_\alpha(n,k)$.
Recall that $I_0:=\min\{i\in\{1,\ldots,n\}\mid S_i/S_{n+1}>t_0\}=nt_0+O_P(n^{1/2})$. Thus
\begin{align*}
 \frac 1{k-1}\sum_{i=1}^{k-1} \Big( \frac{S_i}{S_{n+1}}\vee t_0-t_0\Big)
 & =
\frac 1{k-1}\sum_{i=I_0}^{k-1} \Big( \frac{S_i}{S_{n+1}}\vee t_0-t_0\Big)\\
& = O_P\Big(\frac{k-I_0}k \Big(\frac{S_k}{S_{n+1}}\vee t_0-t_0\Big)\Big)\\
& = o_P\Big(\frac{S_k}{S_{n+1}}\vee t_0-t_0\Big) = o_P(n^{-1/2}).
\end{align*}
Because $H$ vanishes on $(0,t_0]$ and it is differentiable on a right neighborhood of $t_0$ with derivative tending to $h_0$ as $t\downarrow t_0$, we may conclude
\begin{align*}
  \Delta_\alpha(n,k) & = \frac\alpha{k-1} \sum_{i=1}^{k-1} \log\Big( 1+ \frac{H(S_i/S_{n+1})-H(S_k/S_{n+1})}{1+H(S_k/S_{n+1})} \Big)\\
  & = \frac\alpha{k-1} \sum_{i=1}^{k-1} \log\Big(1+ h_0\Big(\frac{S_i}{S_{n+1}}\vee t_0-\frac{S_k}{S_{n+1}}\vee t_0\Big)(1+o(1))\Big)\\
  & = - \frac{\alpha h_0}{k-1} \sum_{i=1}^{k-1} \Big( \frac{S_k}{S_{n+1}}\vee t_0 -\frac{S_i}{S_{n+1}}\vee t_0\Big)(1+o(1))\\
  &  = -\alpha h_0 \Big(\frac{S_k}{S_{n+1}}\vee t_0-t_0\Big)(1+o(1))
\end{align*}
and
\begin{equation}\label{eq:alphaapprox}
\frac{\tilde\alpha_{n,k}-\tilde\alpha^P_k}\alpha = \frac{\alpha/\tilde\alpha_k^P-\alpha/\tilde\alpha_{n,k}}{\alpha^2/(\tilde\alpha_k^P\tilde\alpha_{n,k})} = -\Delta_\alpha(n,k)(1+o_P(1)) = \alpha h_0\Big(\frac{S_k}{S_{n+1}}-t_0\Big)^+(1+o_P(1)).
\end{equation}

Recall that $t_{n,1}\downarrow t_0$ satisfies
 $n^{1/2}
\inf_{t\in[t_{n,1},1]}\tilde D_{\ceil{nt}}\to\infty$. For
$k\in\{\ceil{n(t_0+\eps_n)}, \ldots, $ $\ceil{nt_{n,1}}\}$ let $j_k :=
(k+nt_0)/2$ so that $\log(j_k/k)\to 0$ uniformly. By a Taylor
expansion, we obtain for the second factor in \eqref{eq:Dktrep}
\begin{equation} \label{eq:factor2}
 \Big(\frac{S_{j}}{S_k}\Big)^{(\tilde\alpha_{n,k}-\tilde\alpha^P_k)/\alpha}- 1
  =\alpha h_0\Big(\frac{S_k}{S_{n+1}}-t_0\Big)^+\log\frac{S_{j}}{S_k}(1+o_P(1))
 \end{equation}
if the right hand side tends to 0 in probability. In particular, for $j=j_k$ this factor is stochastically of smaller order than $(S_k/S_{n+1}-t_0)^+$.

In contrast, the third factor has the following asymptotic behavior:
\begin{align} \label{eq:factor3}
  \Big(\frac{1+H(S_j/S_{n+1})}{1+H(S_k/S_{n+1})}\Big)^{-\tilde\alpha_{n,k}}-1
  & = \Big(1+\frac{H(S_j/S_{n+1})-H(S_k/S_{n+1})}{1+H(S_k/S_{n+1})}\Big)^{-\alpha+o_P(1)}-1\nonumber\\
  & = \alpha h_0 \Big(\frac{S_k}{S_{n+1}}\vee t_0-\frac{S_j}{S_{n+1}}\vee t_0\Big)(1+o_P(1)),
\end{align}
which for $j=j_k$ equals
$$ \alpha h_0 \frac{k-j_k}n(1+o_P(1)) = \alpha h_0 \frac 12\Big( \frac kn-t_0\Big)(1+o_P(1)). $$
Combining this with the asymptotic behavior of the other factors, we see that
\begin{align*}
 \limsup_{n\to \infty} & \inf_{k\in \{\ceil{n(t_0+\eps_n)},\ldots,t_{n,1}\}} n^{1/2} \tilde D_k\\
 & \ge
\inf_{k\in \{\ceil{n(t_0+\eps_n)},\ldots,t_{n,1}\}} n^{1/2}\alpha h_0 \frac 12 \Big( \frac kn-t_0\Big)(1+o_P(1))
  \to \infty
\end{align*}
in probability. Together with \eqref{eq:largeDktapprox}, this proves assertion (ii).

It remains to analyze the asymptotic behavior of $\tilde D_k$ for $k=nt_0+n^{1/2}u$ with $u\in[-C,C]$.

Recall from the proof of Theorem \ref{th:Parmain} that, uniformly for $j\in\{1,\ldots,k\}$,
\begin{align}  \label{eq:apprfactor1}
  n^{1/2} \Big(\Big( \frac{S_j}{S_k}\Big)^{\tilde\alpha_k^P/\alpha}-\frac jk\Big) = &
  \frac{W_n(j/n)}{j/n} - \frac{W_n(k/n)}{k/n} +\nonumber\\
   & \hspace*{1cm} +\frac 1{k-1} \sum_{i=1}^{k-1} \Big(\frac{W_n(i/n)}{i/n} - \frac{W_n(k/n)}{k/n} \Big) \log\frac jk +o_P(1)
\end{align}
which tends to 0 uniformly for $j\in\{I_0,\ldots,k\}$ with $I_0:=\min\{i\mid S_i/S_{n+1}>t_0\}$.

Because $(S_k/S_{n+1}-t_0)^+=O_P(n^{-1/2})$, \eqref{eq:factor2} and \eqref{eq:Wncont1} show that
\begin{align} \label{eq:apprfactor2}
  &n^{1/2}  \Big(\Big( \frac{S_j}{S_k}\Big)^{(\tilde\alpha_{n,k}-\tilde\alpha_k^P)/\alpha}-1\Big) \nonumber\\
  & = \alpha h_0 n^{1/2} \Big( \frac{S_k}{S_{n+1}}-t_0\Big)^+\log\frac jk +o_P(1)+O_P(j^{-1/2})\nonumber\\
  & = \alpha h_0 n^{1/2} \Big( \frac{k-nt_0+W(k)-t_0W(n)+O(\log n)}{n+W(n)+O(\log n)}\Big)^+ \log\frac jk (1+o_P(1))+O_P(j^{-1/2})\nonumber\\
  & = \alpha h_0 \big( u+W_n(t_0)-t_0W_n(1)\big)^+ \log\frac jk(1+ o_P(1))+O_P(j^{-1/2}).
\end{align}
This expression converges to 0 uniformly for $j\in\{I_0,\ldots,k\}$, too.

Likewise, in view of \eqref{eq:factor3}, we have
\begin{equation} \label{eq:apprfactor3}
   \Big(\frac{1+H(S_j/S_{n+1})}{1+H(S_k/S_{n+1})}\Big)^{-\tilde\alpha_{n,k}}-1 = \alpha h_0\Big(\frac{S_k}{S_{n+1}}\vee t_0-\frac{S_j}{S_{n+1}}\vee t_0\Big)+ o_P(n^{-1/2}).
\end{equation}
Here, $S_k/S_{n+1}\vee t_0-S_j/S_{n+1}\vee t_0$ equals $(S_k/S_{n+1}-t_0)^+=n^{-1/2}(u+W_n(t_0)-t_0W_n(1)+o_P(1))^+$ for $j<I_0$, and it is decreasing  in $j\in\{I_0,\ldots,k\}$ with value 0 for $j=k$.

Combining \eqref{eq:apprfactor1}-- \eqref{eq:apprfactor3} with \eqref{eq:Dktrep}, we arrive at
\begin{align*}
 n^{1/2}& \tilde D_k \nonumber\\
  = & \max_{1\le j\le k}n^{1/2}\bigg|\bigg(\frac jk + n^{-1/2}\Big(\frac{W_n(j/n)}{j/n} - \frac{W_n(k/n)}{k/n} + \\
  & \hspace*{3cm} +\frac 1{k-1} \sum_{i=1}^{k-1} \Big(\frac{W_n(i/n)}{i/n} - \frac{W_n(k/n)}{k/n} \Big) \log\frac jk +o_P(1)\Big)\bigg)\times\\
 & \times \bigg(1+n^{-1/2}\Big(\alpha h_0 \big( u+W_n(t_0)-t_0W_n(1)\big)^+ \log\frac jk +o_P(1)+O_P(j^{-1/2})\Big)\bigg)\times\\
 & \times  \bigg(1+n^{-1/2}\Big(\alpha h_0 \big( u+W_n(t_0)-t_0W_n(1)\big)^+ +o_P(1)\Big)\bigg)-\frac jk\bigg|\\
  = & \max_{1\le j\le I_0} \bigg| \frac{W_n(j/n)}{j/n} - \frac{W_n(k/n)}{k/n} + \frac 1{k-1} \sum_{i=1}^{k-1} \Big(\frac{W_n(i/n)}{i/n} - \frac{W_n(k/n)}{k/n} \Big) \log\frac jk \\
 & \hspace*{1cm}+ \Big(\alpha h_0 \big( u+W_n(t_0)-t_0W_n(1)\big)^+\frac jk\Big(1+\log\frac jk\Big)\bigg| + o_P(1)
\end{align*}
uniformly for $u\in[-C,C]$. Now assertion (iii) follows by the arguments given in the last part of the proof of Theorem \ref{th:Parmain}.

\end{proofof}

\section*{Acknowledgments}
Sidney Resnick and Tiandong Wang were partly supported by
US Army MURI grant W911NF-12-1-0385 to Cornell University. Holger Drees was partly supported by DFG grant DR 271/6-2 as part of the Research Unit 1735. The authors appreciate the thoughtful, helpful comments from referees and editors.

\bigskip

\bibliographystyle{siamplain}
\bibliography{MDSP}
\end{document}